# RECURRENCE OF EDGE-REINFORCED RANDOM WALK ON A TWO-DIMENSIONAL GRAPH

By Franz Merkl and Silke W. W. Rolles

We consider a linearly edge-reinforced random walk on a class of two-dimensional graphs with constant initial weights. The graphs are obtained from $\mathbb{Z}^2$ by replacing every edge by a sufficiently large, but fixed number of edges in series. We prove that the linearly edge-reinforced random walk on these graphs is recurrent. Furthermore, we derive bounds for the probability that the edge-reinforced random walk hits the boundary of a large box before returning to its starting point.

**1. Introduction.** In his paper [10], Pemantle studies recurrence and transience of a linearly edge-reinforced random walk on trees with a parameter $\Delta > 0$. At the end of this paper, he writes: "Questions of reinforced random walk on other graphs are still wide open. Diaconis originally asked me about the $d$-dimensional integer lattice $\mathbb{Z}^d$. I believe it is not even known whether there is a $\Delta > 0$ for which the reinforced random walk on $\mathbb{Z}^2$ is recurrent!"

This turned out to be a hard problem which is unresolved. However, in this article, we solve a variant of this problem. For a class of fully two-dimensional translationally symmetric graphs and sufficiently small constant initial weights, we show that the edge-reinforced random walk visits every vertex infinitely often with probability one. This is the first time that recurrence of the linearly edge-reinforced random walk is proven for a fully two-dimensional graph.

*Comparison with previous work.* Earlier versions of the technique presented here have been used in [7] to prove recurrence for the edge-reinforced random walk on a large class of one-dimensional reflection-symmetric periodic graphs, including, among others, $\mathbb{Z} \times \{1, \ldots, d\}$ for $d \in \mathbb{N}$. The technique can be also used to prove some bounds for the random environment corresponding to the edge-reinforced random walk on $\mathbb{Z}^2$; see [9].

Given a graph $G$, its $r$-diluted version is the graph obtained from $G$ by replacing every edge by $r$ edges in series. The basic new idea of the present paper is that this dilution of a graph makes recurrence simpler to prove: Already a simplification and generalization of the techniques from [9] to diluted two-dimensional graphs suffices to prove recurrence.

For a more detailed comparison of previous papers with the present work, we refer the reader to the concluding remarks in Section 10. An overview of earlier work on the linearly edge-reinforced random walk can be found in [6]. For a recent overview of general processes with reinforcement, we refer the reader to [11].

*Definition of linearly edge-reinforced random walk.* Given $r \in \mathbb{N}$, we consider the graph $G_r = (V_r, E_r)$ with vertex set

$$(1.1) \qquad V_r = \{(x_1, x_2) \in \mathbb{Z}^2 : x_1 \in r\mathbb{Z} \text{ or } x_2 \in r\mathbb{Z}\}$$

and edge set

$$(1.2) \qquad E_r = \{\{u, v\} \subset V_r : |u - v| = 1\}.$$

Here $|x|$ denotes the Euclidian norm of $x$. Note that the edges are *undirected*. A finite piece of $G_r$ for $r = 4$ is shown in Figure 1.

Let $\mathbf{0} := (0, 0)$. The *linearly edge-reinforced random walk* (ERRW) on $G_r$ with constant initial weights $a > 0$ and starting point $\mathbf{0}$ is a stochastic process $(X_t)_{t \in \mathbb{N}_0}$ with law $P_{\mathbf{0},a}^{G_r}$ defined as follows: At every discrete time $t \in \mathbb{N}_0$, every edge $e \in E_r$ is assigned a strictly positive number $w_t(e)$ as a weight. Initially, all weights are equal to $a$:

$$(1.3) \qquad w_0(e) = a \qquad \text{for all } e \in E_r.$$

The edge-reinforced random walker starts in the vertex $\mathbf{0}$ at time 0:

$$(1.4) \qquad P_{\mathbf{0},a}^{G_r}[X_0 = \mathbf{0}] = 1.$$

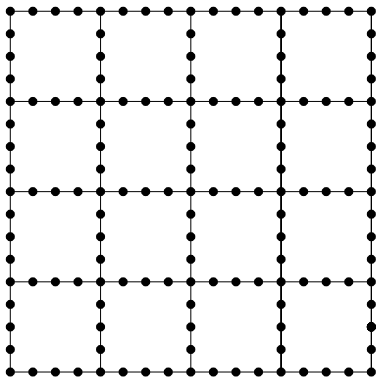

FIG. 1. *A finite piece of the graph $G_4$.*

At each discrete time $t \in \mathbb{N}_0$, the random walker jumps randomly from its current position $X_t$ to a neighboring vertex $v$ with probability proportional to the current weight of the connecting edge $\{X_t, v\}$:

$$(1.5) \quad P_{\mathbf{0},a}^{G_r}[X_{t+1} = v \mid X_0, X_1, \ldots, X_t] = \frac{w_t(\{X_t, v\})}{\sum_{e \in E_r : e \ni X_t} w_t(e)} 1_{\{\{X_t, v\} \in E_r\}}.$$

The weight of the traversed edge is immediately increased by 1, whereas all other weights remain unchanged:

$$(1.6) \quad w_{t+1}(e) = w_t(e) + 1_{\{\{X_t, X_{t+1}\} = e\}} \quad \text{for all } e \in E_r.$$

Thus, the weight of edge $e$ at time $t$ equals the initial weight increased by the number of times the reinforced random walker has traversed $e$ up to time $t$:

$$(1.7) \quad w_t(e) = a + \sum_{s=0}^{t-1} 1_{\{\{X_s, X_{s+1}\} = e\}}.$$

We call a path $(v_0, v_1, v_2, \ldots)$ in a graph *admissible* if $\{v_{t-1}, v_t\}$ is an edge in the graph for all $t \in \mathbb{N}$. We realize $P_{\mathbf{0},a}^{G_r}$ as a probability measure on the set $\Sigma \subseteq V_r^{\mathbb{N}_0}$ of admissible paths in $G_r$, not necessarily starting in $\mathbf{0}$, endowed with the $\sigma$-field generated by the canonical projections $X_t : \Sigma \to V_r$, $t \in \mathbb{N}_0$.

*Main results.* We prove the following:

THEOREM 1.1 (Recurrence). *For all $r \in \mathbb{N}$ with $r \geq 130$ and all $a \in (0, (r-129)/256)$, the linearly edge-reinforced random walk on $G_r$ with constant initial weights $w_0 \equiv a$ visits all vertices infinitely often with probability one.*

Note that small initial weights $a$ correspond to strong reinforcement. This can be seen more intuitively by taking the initial weights equal to one, but the weight increment at each time step equal to $1/a$, which defines the same law of the random walk. In this sense, the recurrence result applies to strong linear reinforcement.

In order to prove recurrence, we derive bounds for the probability that the edge-reinforced random walk hits the boundary of a large box before returning to its starting point. Let us introduce some notation before we state the result: For $A \subseteq V_r$, let

$$\tau_A := \inf\{t \geq 1 : X_t \in A\}$$

be the hitting time of $A$. If $A = \{v\}$ contains just one vertex, we simply write $\tau_v$ instead of $\tau_{\{v\}}$. Let

$$(1.8) \quad L_r := r\mathbb{Z}^2$$

be the set of "four-way-crossings" in the graph $G_r$. For $(v_1, v_2) \in V_r$, set $|(v_1, v_2)|_\infty := \max\{|v_1|, |v_2|\}$.

We prove the following:

THEOREM 1.2 (Hitting probabilities). *For all $r \in \mathbb{N}$ with $r \geq 130$ and all initial weights $a \in (0, (r-129)/256)$, there exist constants $l_0 = l_0(r,a) \in \mathbb{N}$ and $\xi = \xi(r,a) > 0$, such that the following hold:*

(a) *For all $\ell \in L_r$ with $|\ell|_\infty \geq l_0$, the probability that the edge-reinforced random walker hits $\ell$ before returning to its starting point satisfies*

$$(1.9) \qquad P_{\mathbf{0},a}^{G_r}[\tau_\ell < \tau_{\mathbf{0}}] \leq \left(\frac{r}{|\ell|_\infty}\right)^{1+\xi}.$$

(b) *As a consequence, for all $l \geq l_0$, the probability that the edge-reinforced random walker hits a vertex in the set $\mathcal{V}_l := \{v \in V_r : |v|_\infty = rl\}$ (i.e. a vertex in the intersection of $V_r$ with the boundary of the box $[-rl, rl]^2$) before returning to $\mathbf{0}$ can be bounded as follows:*

$$(1.10) \qquad P_{\mathbf{0},a}^{G_r}[\tau_{\mathcal{V}_l} < \tau_{\mathbf{0}}] \leq 8l^{-\xi}.$$

*Notation.* For $a, b \in \mathbb{R}$, set $a \wedge b = \min\{a, b\}$ and $a \vee b = \max\{a, b\}$. Let $\lfloor a \rfloor$ denote the largest integer $\leq a$, and let $\lceil a \rceil$ denote the smallest integer $\geq a$. The cardinality of a set $A$ is denoted by $|A|$. For two probability measures $P$ and $Q$ on the same space, $P \ll Q$ denotes that $P$ is absolutely continuous with respect to $Q$.

**2. Main lemmas and global structure of the paper.** All the theorems in Section 1 rely on the bound (1.9) of the hitting probabilities. The main work in this paper is required to prove this bound.

*Finite boxes with periodic boundary conditions.* The bound (1.9) is based on an analysis of the edge-reinforced random walk on *finite* boxes with periodic boundary conditions, in the limit as the box size goes to infinity. We now introduce these finite boxes and some related notation.

Given $r \in \mathbb{N}$ and $i \in \mathbb{N}$ with $i > 1$, we define a finite variant $G_r^{(i)} = (V_r^{(i)}, E_r^{(i)})$ of $G_r$ with periodic boundary conditions as follows: Set

$$(2.1) \qquad V_r^{(i)} := \{v + ri\mathbb{Z}^2 : v \in V_r\} \subseteq \mathbb{Z}^2/(ri\mathbb{Z}^2),$$

$$(2.2) \qquad E_r^{(i)} := \{\{u + ri\mathbb{Z}^2, v + ri\mathbb{Z}^2\} : \{u, v\} \in E_r\}.$$

Whenever there is no risk of confusion, we identify $V_r^{(i)}$ with its set of representatives

(2.3) $\tilde{V}_r^{(i)} := \{(v_1, v_2) \in V_r : r\lfloor i/2 \rfloor - ri < v_j \leq r\lfloor i/2 \rfloor \text{ for } j = 1, 2\} \subset V_r$.

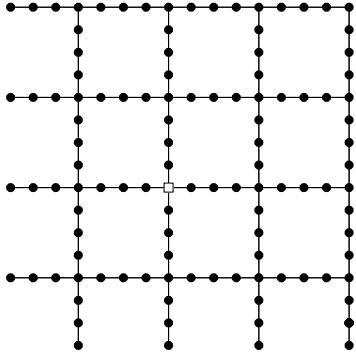

FIG. 2. *The graph $(\tilde{V}_4^{(4)}, \tilde{E}_4^{(4)})$. The white square marks the origin.*

Furthermore, we identify the subset

(2.4) $\quad \hat{E}_r^{(i)} := \{\{u + ri\mathbb{Z}^2, v + ri\mathbb{Z}^2\} : \{u,v\} \in E_r, u, v \in \tilde{V}_r^{(i)}\} \subset E_r^{(i)}$

with its set of representatives

(2.5) $\qquad \tilde{E}_r^{(i)} := \{\{u,v\} \in E_r : u, v \in \tilde{V}_r^{(i)}\} \subset E_r.$

Note that the inclusion in (2.4) is strict. We call all edges in $E_r^{(i)} \setminus \hat{E}_r^{(i)}$ *periodically closing* edges. For instance, the edge $\{(r\lfloor i/2 \rfloor, 0) + ri\mathbb{Z}^2, (r\lfloor i/2 \rfloor - ri + 1, 0) + ri\mathbb{Z}^2\}$ is periodically closing as it is contained in $E_r^{(i)}$, but not in $\hat{E}_r^{(i)}$.

Figure 2 illustrates the graph $(\tilde{V}_r^{(i)}, \tilde{E}_r^{(i)})$ for $r = i = 4$.

Finally, we set

(2.6) $\qquad L_r^{(i)} := \{v + ri\mathbb{Z}^2 : v \in L_r\} \subseteq V_r^{(i)}.$

Thus, $L_r^{(i)}$ consists of the vertices in $G_r^{(i)}$ of degree 4.

To prove our main results, we analyze the edge-reinforced random walk on the finite boxes $G_r^{(i)}$ with a constant initial weight $a > 0$.

*Linearly edge-reinforced on finite graphs.* Since most parts of our arguments are quite robust, we work in this section with the edge-reinforced random walk on a general finite connected undirected graph $G = (V, E)$ with possibly edge-dependent initial weights $a_e > 0$, $e \in E$. In Section 8, we specialize to the case $G = G_r^{(i)}$.

The edge-reinforced random walk on a finite graph $G = (V, E)$ is defined in analogy to the definition on $G_r$ given in the Introduction. Given a starting vertex $v_0 \in V$ and initial weights $\mathbf{a} = (a_e)_{e \in E}$, it is a stochastic process

$(X_t)_{t\in\mathbb{N}_0}$ taking values in $V$ with law $P^G_{v_0,\mathbf{a}}$, specified by

$$w_t(e) = a_e + \sum_{s=0}^{t-1} 1_{\{\{X_s,X_{s+1}\}=e\}} \qquad (e \in E, t \in \mathbb{N}_0), \tag{2.7}$$

$$P^G_{v_0,\mathbf{a}}[X_0 = v_0] = 1, \tag{2.8}$$

$$P^G_{v_0,\mathbf{a}}[X_{t+1} = v \mid X_0, X_1, \ldots, X_t] = \frac{w_t(\{X_t,v\})}{\sum_{e\in E: e \ni X_t} w_t(e)} 1_{\{\{X_t,v\}\in E\}} \tag{2.9}$$

$$(v \in V, t \in \mathbb{N}_0).$$

Again, we realize $P^G_{v_0,\mathbf{a}}$ as a probability measure on the set $\Sigma_G \subseteq V^{\mathbb{N}_0}$ of admissible paths in $G$, not necessarily starting in $v_0$, endowed with the $\sigma$-field generated by the canonical projections $X_t : \Sigma_G \to V$, $t \in \mathbb{N}_0$.

*ERRW as a mixture.* It is known that this edge-reinforced random walk on the finite graph $G$ is a mixture of reversible Markov chains with a unique mixing measure. Let us review this result. Let $\Omega = (0,\infty)^E$ be the space of strictly positive edge weights on $G$. For $x = (x_e)_{e \in E} \in \Omega$ and $v \in V$, define

$$x_v := \sum_{\substack{e \in E: \\ e \ni v}} x_e. \tag{2.10}$$

For $v_0 \in V$, we denote by $Q_{v_0,x}$ the distribution of the Markovian random walk on $G$ which starts in $v_0$ and jumps from $u$ to a neighboring point $v$ with probability proportional to the weight $x_e$ of the connecting edge $e = \{u,v\}$:

$$Q_{v_0,x}[X_0 = v_0] = 1, \tag{2.11}$$

$$Q_{v_0,x}[X_{t+1} = v \mid X_0, X_1, \ldots, X_t] = \frac{x_{\{X_t,v\}}}{x_{X_t}} 1_{\{\{X_t,v\}\in E\}} \tag{2.12}$$

$$(v \in V, t \in \mathbb{N}_0).$$

Fix a reference edge $e_0 \in E$. At this point, the edge $e_0$ is arbitrary, but from Section 5 onward, it will be convenient to assume $v_0 \in e_0$. We introduce the set of weights

$$\Omega_{e_0} := \{x \in \Omega : x_{e_0} = 1\} \tag{2.13}$$

having the weight of the reference edge $e_0$ normalized.

LEMMA 2.1 (Mixture of reversible Markov chains, Corollary 3.1 of [12]). *For all $v_0 \in V$, $\mathbf{a} \in (0,\infty)^E$, and $e_0 \in E$, there exists a unique probability measure $\mathbb{Q}^G_{v_0,\mathbf{a},e_0}$ on $\Omega_{e_0}$ with the following properties: The distribution $P^G_{v_0,\mathbf{a}}$* of the edge-reinforced random walk on $G$ starting in $v_0$ with initial weights $\mathbf{a}$ satisfies

$$P^G_{v_0,\mathbf{a}}[A] = \int_{\Omega_{e_0}} Q_{v_0,x}[A] \mathbb{Q}^G_{v_0,\mathbf{a},e_0}(dx) \tag{2.14}$$

for all events $A \subseteq \Sigma_G$ of admissible paths in $G$.

We call $\mathbb{Q}^G_{v_0,\mathbf{a},e_0}$ the "*mixing measure.*"

The proof of this lemma in [12] is based on a de Finetti theorem for *reversible* Markov chains, which relies on a de Finetti theorem for Markov chains by Diaconis and Freedman [2].

*The main lemma—concrete version.* All theorems of this paper are based on a lemma, which we phrase in this section and prove in the subsequent sections. There are two versions of this main lemma: a "concrete" version, which works with finite boxes with periodic boundary conditions and constant initial weights, and an "abstract" version working with general finite graphs and general initial weights having a reflection symmetry. The concrete version is a consequence of the abstract version.

For the edge-reinforced random walk on $G_r^{(i)}$ with constant initial weights $\mathbf{a} = (a_e \equiv a)_{e \in E_r^{(i)}}$, starting point $\mathbf{0}$, and reference edge $e_0 = \{\mathbf{0}, (0,1)\}$, we denote the mixing measure from Lemma 2.1 by

$$\mathbb{Q}_\mathbf{0}^{(i)} := \mathbb{Q}^{G_r^{(i)}}_{\mathbf{0},\mathbf{a},e_0}. \tag{2.15}$$

The initial weight $a$ is omitted from the notation because it is fixed throughout the paper.

Roughly speaking, the concrete version of the main lemma bounds in a stochastic sense the ratio of two vertex weights $x_\ell$ and $x_\mathbf{0}$ with respect to the mixing measure $\mathbb{Q}_\mathbf{0}^{(i)}$. For $\ell \in L_r$, let $i_0(\ell) \in \mathbb{N}$ be large enough that, for all $i \geq i_0(\ell)$, the vertex $\ell$ belongs to $\tilde{V}_r^{(i)}$ and has vertex degree 4 in the graph $(\tilde{V}_r^{(i)}, \tilde{E}_r^{(i)})$.

LEMMA 2.2 (Main lemma—concrete version). *Let $r \in \mathbb{N}$ with $r \geq 130$ and let the initial weight satisfy $a \in (0, (r-129)/256)$. There exist constants $l_0 = l_0(r,a) \in \mathbb{N}$ and $\xi = \xi(r,a) > 0$, such that, for all $\ell \in L_r$ with $|\ell|_\infty \geq l_0$ and $i \geq i_0(\ell)$, the following estimate holds:*

$$E_{\mathbb{Q}_\mathbf{0}^{(i)}}\left[\left(\frac{x_\ell}{x_\mathbf{0}}\right)^{1/4}\right] \leq \left(\frac{r}{|\ell|_\infty}\right)^{1+\xi}. \tag{2.16}$$

*The main lemma—abstract version.* The abstract version of the main lemma works with a general finite connected graph $G = (V, E)$ with arbitrary initial weights $\mathbf{a} \in (0, \infty)^E$. The weighted graph $(G, \mathbf{a})$ needs to have an additional *reflection symmetry*, which we introduce now.

We call a bijection $f : V \to V$ on the vertices an *automorphism* of the *weighted* graph $(G, \mathbf{a})$, if the following hold:

- For all $u, v \in V$, one has $\{f(u), f(v)\} \in E$ if and only if $\{u, v\} \in E$.
- The edge weights are preserved: $a_{\{f(u), f(v)\}} = a_{\{u,v\}}$ holds for all $\{u, v\} \in E$.

We extend $f$ to the set of edges as follows: $f(e) = \{f(u), f(v)\}$ for $e = \{u, v\} \in E$.

In the abstract version of the main lemma, the vertices $\mathbf{0}$ and $\ell$ from the concrete version are replaced by two vertices $v_0, v_1 \in V$ with the following symmetry property:

ASSUMPTION 2.3. *We take two vertices $v_0, v_1 \in V$, a reference edge $e_0 \in E$, and a function $\varphi : E \to [0, 1]$, and we assume the following:*

(a) *The vertices $v_0, v_1 \in V$ satisfy $v_0 \neq v_1$, and $v_0$ is not adjacent to $v_1$, that is, $\{v_0, v_1\} \notin E$.*
(b) *The reference edge $e_0$ is incident to $v_0$: $v_0 \in e_0$.*
(c) *There exists an automorphism $f$ of the weighted graph $(G, \mathbf{a})$ with $f(v_0) = v_1$ and $f(v_1) = v_0$.*
(d) *We have $\varphi(e) = 0$ for all $e \in E$ with $v_0 \in e$, and $\varphi(e) = 1$ for all $e \in E$ with $v_1 \in e$.*

REMARK 2.4. In the special case of the periodic boxes $G_r^{(i)}$, endowed with constant initial weights $a$, parts (a)–(c) of this assumption are satisfied for $r > 1$, $v_0 = \mathbf{0}$, $v_1 = \ell \in L_r^{(i)} \setminus \{\mathbf{0}\}$, and $e_0 = \{\mathbf{0}, (0, 1)\}$. In particular, the two vertices $\mathbf{0}$ and $\ell$ have both four neighbors; recall Definition (2.6) of $L_r^{(i)}$. In this special case, we may take $f$ as the reflection

$$(2.17) \qquad f : V \to V, \qquad v \mapsto \ell - v.$$

For $v \in V$, set

$$(2.18) \qquad a_v := \sum_{e \in E : e \ni v} a_e.$$

LEMMA 2.5 (Main lemma—abstract version). *Under Assumption 2.3, one has the bound*

$$(2.19) \qquad E_{\mathbb{Q}_{v_0, \mathbf{a}, e_0}^G}\left[\left(\frac{x_{v_1}}{x_{v_0}}\right)^{1/4}\right] \leq \exp\left(-\frac{1}{32 S_\varphi}\right),$$

*where*

$$(2.20) \qquad S_\varphi = \sum_{v \in V} \frac{a_v + 1}{2} \max_{\substack{e, e' \in E: \\ v \in e, v \in e'}} (\varphi(e) - \varphi(e'))^2.$$

*This article is organized as follows.* In Section 3, we review an explicit description of the mixing measure $\mathbb{Q}^G_{v_0, \mathbf{a}, e_0}$ from Lemma 2.1. This description was discovered by Coppersmith and Diaconis [1] in 1986. Later, Keane and Rolles [4] published a proof of their formula. Section 3 also contains some properties of the mixing measure, in particular, scaling properties, which are needed in the present paper.

In Section 4, we interpolate between the two random environment distributions corresponding to two different starting points. We also show that properties of the interpolations are inherited from the scaling properties of the mixing measure. In particular, normalizing constants remain unchanged after scaling.

Proving the bounds for the weights boils down to proving an upper bound for these normalizing constants; see (4.13) below. The upper bound for the normalizing constant is then obtained via a general variational principle, Lemma 5.1. This variational principle is proved in Section 5. It relies only on the Assumptions 2.3 for the underlying weighted graph.

The variational principle is then applied to a deformation of the interpolated random environments. In Sections 5–7, this leads to a proof of Lemma 2.5, the main lemma in the context of a general finite connected graph.

In Section 8, we verify that the assumptions of the abstract version of the main lemma (Lemma 2.5) are satisfied for the edge-reinforced random walk on the boxes $G_r^{(i)}$ with constant initial weights $a$, provided $r$ and $a$ are chosen appropriately. In this way, we deduce Lemma 2.2 from Lemma 2.5.

Finally, in Section 9, we derive the estimates for the hitting probabilities of the edge-reinforced random walk (Theorem 1.2) from Lemma 2.5. Recurrence is shown using these estimates for the hitting probabilities.

*The following results from other papers are used.* The present paper uses only a few results from our previous articles. One essential ingredient is Lemma 3.4, which is a variant of the formula by Coppersmith and Diaconis [1]. The version of this formula used in the present paper is taken from the Appendix of [9]. Besides this description of the mixing measure, we cite only a general bound for the tails of ratios of two edge weights (Theorem 2.4 in [8]) and the fact that the edge-reinforced random walk visits every vertex infinitely often with probability one if and only if it returns with probability one to its starting point (Theorem 2.1 of [8]).

**3. Properties of the mixing measure.** Throughout Sections 3–7, we consider the edge-reinforced random walk on a finite connected graph $G$ with initial weights $\mathbf{a} \in (0,\infty)^E$. Since the underlying graph $G$ is fixed, we suppress the dependence on $G$ in the notation. In particular, we introduce the following abbreviation for the mixing measure $\mathbb{Q}^G_{v_0,\mathbf{a},e_0}$ from Lemma 2.1:

$$(3.1) \qquad \mathbb{Q}_{v_0,\mathbf{a},e_0} := \mathbb{Q}^G_{v_0,\mathbf{a},e_0}.$$

The mixing measure $\mathbb{Q}_{v_0,\mathbf{a},e_0}$ has an explicit description, discovered by Coppersmith and Diaconis [1, 3]. We review it in this section. This description reveals a nice scaling property of the mixing measure. In order to state it, we need some preparations.

Denote by $\delta_1$ the Dirac measure at 1. We endow $\Omega_{e_0}$ with the reference measure

$$(3.2) \quad \rho_{e_0}(dx) := \delta_1(dx_{e_0}) \prod_{e \in E \setminus \{e_0\}} \frac{dx_e}{x_e} = \delta_1(dx_{e_0}) \prod_{e \in E \setminus \{e_0\}} d(\log x_e).$$

The measure $\rho_{e_0}$ is just the Lebesgue measure on a logarithmic scale. More precisely, the image of $\rho_{e_0}$ under the bijection

$$(3.3) \qquad L_{e_0} : \Omega_{e_0} \to \mathbb{R}^{E \setminus \{e_0\}}, \qquad (x_e)_{e \in E} \mapsto (\log x_e)_{e \in E \setminus \{e_0\}}$$

is the Lebesgue measure on $\mathbb{R}^{E \setminus \{e_0\}}$. Given a different reference edge $e_1 \in E$, we introduce the map

$$(3.4) \qquad g_{e_0,e_1} : \Omega_{e_0} \to \Omega_{e_1}, \qquad g_{e_0,e_1}(x) = (x_e/x_{e_1})_{e \in E}.$$

It serves to change the normalization. The map $g_{e_0,e_1}$ is a bijection with inverse $g_{e_1,e_0}$.

LEMMA 3.1. *For all $e_0, e_1 \in E$, $e_0 \neq e_1$, and all Borel sets $A \subseteq \Omega_{e_1}$, we have*

$$(3.5) \qquad \rho_{e_0}[\{x \in \Omega_{e_0} : g_{e_0,e_1}(x) \in A\}] = \rho_{e_1}[A].$$

*In words, $\rho_{e_1}$ is the image measure of $\rho_{e_0}$ with respect to $g_{e_0,e_1}$.*

PROOF. It suffices to prove the lemma for sets of the form $A = \prod_{e \in E} A_e$ with $A_{e_1} = \{1\}$ and $A_e = (b_e, c_e]$ for $e \in E \setminus \{e_1\}$, where $0 < b_e < c_e$. Then, for $x \in \Omega_{e_0}$, we have $g_{e_0,e_1}(x) \in A$ if and only if the following hold: $1/c_{e_0} \leq x_{e_1} < 1/b_{e_0}$, and $x_{e_1} b_e < x_e \leq x_{e_1} c_e$ for $e \in E \setminus \{e_0, e_1\}$. Thus, the left-hand side in (3.5) equals

$$(3.6) \quad \int_{1/c_{e_0}}^{1/b_{e_0}} \left( \prod_{e \in E \setminus \{e_0,e_1\}} \int_{x_{e_1} b_e}^{x_{e_1} c_e} \frac{dx_e}{x_e} \right) \frac{dx_{e_1}}{x_{e_1}} = \int_{1/c_{e_0}}^{1/b_{e_0}} \left( \prod_{e \in E \setminus \{e_0,e_1\}} \log \frac{c_e}{b_e} \right) \frac{dx_{e_1}}{x_{e_1}}$$

$$= \prod_{e \in E \setminus \{e_1\}} \log \frac{c_e}{b_e} = \rho_{e_1}[A].$$

ERRW ON A 2-DIMENSIONAL GRAPH 11□

The following function $\Phi_{v_0,\mathbf{a}}$ plays an essential role in the description of the mixing measure.

DEFINITION 3.2. A spanning tree $T$ of $G$ is a maximal connected subgraph of $G$ without cycles, viewed as a set $T \subseteq E$ of edges. Let $\mathcal{T}$ denote the set of spanning trees of the graph $G$. We introduce the map $\Phi_{v_0,\mathbf{a}} : \Omega \to (0,\infty)$, given by

$$(3.7) \qquad \Phi_{v_0,\mathbf{a}}((x_e)_{e\in E}) = \frac{\prod_{e\in E} x_e^{a_e}}{x_{v_0}^{a_{v_0}/2} \prod_{v\in V\setminus\{v_0\}} x_v^{(a_v+1)/2}} \sqrt{\sum_{T\in\mathcal{T}} \prod_{e\in T} x_e}.$$

LEMMA 3.3 (Scaling property). *For all $\lambda \in (0,\infty)$ and all $x = (x_e)_{e\in E} \in \Omega$, we have*

$$(3.8) \qquad \Phi_{v_0,\mathbf{a}}(\lambda x) = \Phi_{v_0,\mathbf{a}}(x).$$

*As a consequence, for all $e_0, e_1 \in E$ and all $x \in \Omega_{e_0}$,*

$$(3.9) \qquad \Phi_{v_0,\mathbf{a}}(g_{e_0,e_1}(x)) = \Phi_{v_0,\mathbf{a}}(x).$$

PROOF. Let $\lambda \in (0,\infty)$. Any spanning tree $T \in \mathcal{T}$ has $|V|-1$ edges. Hence, for any $x \in \Omega$, we get

$$(3.10) \qquad \Phi_{v_0,\mathbf{a}}(\lambda x) = \lambda^A \Phi_{v_0,\mathbf{a}}(x)$$

with the exponent

$$(3.11) \quad A = \sum_{e\in E} a_e - \frac{a_{v_0}}{2} - \sum_{v\in V\setminus\{v_0\}} \frac{a_v+1}{2} + \frac{|V|-1}{2} = \sum_{e\in E} a_e - \frac{1}{2}\sum_{v\in V} a_v = 0;$$

for the last equality, we used the definition of $a_v$. This proves (3.8). By the definition (3.4) of $g_{e_0,e_1}$, the identity (3.9) is just the scaling property (3.8) in the special case $\lambda = 1/x_{e_1}$. □

The mixing measure from Lemma 2.1 can be described as follows:

LEMMA 3.4 (Description of the mixing measure, Lemma 9.1 of [9]). *For all $v_0 \in V$, $\mathbf{a} \in (0,\infty)^E$, and $e_0 \in E$, the distribution $\mathbb{Q}_{v_0,\mathbf{a},e_0}$ of the random environment from Lemma 2.1 is absolutely continuous with respect to $\rho_{e_0}$ with the density*

$$(3.12) \qquad \frac{d\mathbb{Q}_{v_0,\mathbf{a},e_0}}{d\rho_{e_0}}(x) = \frac{1}{z_{v_0,e_0}} \Phi_{v_0,\mathbf{a}}(x),$$

*$x \in \Omega_{e_0}$, where $z_{v_0,e_0} > 0$ denotes a normalizing constant.*

Although the distribution $\mathbb{Q}_{v_0,\mathbf{a},e_0}$ *does* depend on the initial weights $\mathbf{a}$, we usually drop the subscript $\mathbf{a}$, writing

$$(3.13) \qquad \mathbb{Q}_{v_0,e_0} = \mathbb{Q}_{v_0,\mathbf{a},e_0}$$

to simplify the notation, since $\mathbf{a}$ is being fixed.

LEMMA 3.5.

(a) *For all $e_0, e_1 \in E$ and all measurable functions $f : \Omega_{e_1} \to [0,\infty]$, one has*

$$(3.14) \qquad \int_{\Omega_{e_0}} f(g_{e_0,e_1}(x)) \Phi_{v_0,\mathbf{a}}(x) \rho_{e_0}(dx) = \int_{\Omega_{e_1}} f(x) \Phi_{v_0,\mathbf{a}}(x) \rho_{e_1}(dx).$$

*In other words, the image of the measure $\Phi_{v_0,\mathbf{a}}\, d\rho_{e_0}$ on $\Omega_{e_0}$ with respect to the map $g_{e_0,e_1} : \Omega_{e_0} \to \Omega_{e_1}$ equals the measure $\Phi_{v_0,\mathbf{a}}\, d\rho_{e_1}$ on $\Omega_{e_1}$.*

(b) *As a consequence, the normalizing constant $z_{v_0,e_0}$ does not depend on the choice of the reference edge:*

$$(3.15) \qquad z_{v_0,e_0} = z_{v_0,e_1}.$$

(c) *Under $\mathbb{Q}_{v_0,e_0}$, the scaling $(x_e/x_{e_1})_{e \in E}$ of the random environment according to $g_{e_0,e_1}$ has the law $\mathbb{Q}_{v_0,e_1}$.*

(d) *The joint laws of the family $(x_e/x_{\tilde{e}}, x_e/x_v, x_v/x_{\tilde{v}})_{e,\tilde{e}\in E, v,\tilde{v}\in V}$ of all ratios of weights with respect to $\mathbb{Q}_{v_0,e_0}$ and with respect to $\mathbb{Q}_{v_0,e_1}$ coincide.*

PROOF. Let $e_0, e_1 \in E$. Using first the scaling property (3.9) and then the transformation formula and Lemma 3.1, we obtain the claim (3.14) in (a) as follows:

$$(3.16) \qquad \begin{aligned} \int_{\Omega_{e_0}} & f(g_{e_0,e_1}(x)) \Phi_{v_0,\mathbf{a}}(x) \rho_{e_0}(dx) \\ &= \int_{\Omega_{e_0}} f(g_{e_0,e_1}(x)) \Phi_{v_0,\mathbf{a}}(g_{e_0,e_1}(x)) \rho_{e_0}(dx) \\ &= \int_{\Omega_{e_1}} f(x) \Phi_{v_0,\mathbf{a}}(x) \rho_{e_1}(dx). \end{aligned}$$

The claim (3.15) in (b) is just the special case $f \equiv 1$ of the first claim (3.14).

Part (c) is an immediate consequence of (a) and (b), since $d\mathbb{Q}_{v_0,e_0} = (z_{v_0,e_0})^{-1} \Phi_{v_0,\mathbf{a}}\, d\rho_{e_0}$ and $d\mathbb{Q}_{v_0,e_1} = (z_{v_0,e_1})^{-1} \Phi_{v_0,\mathbf{a}}\, d\rho_{e_1}$.

Part (d) is a consequence of part (c). To see this, take $x \in \Omega_{e_0}$ and $y = g_{e_0,e_1}(x) \in \Omega_{e_1}$. Using the notation $y_v = \sum_{e \in E : e \ni v} y_e$ for $v \in V$ in analogy to $x_v = \sum_{e \in E : e \ni v} x_e$, we get the claim (d) from

$$(3.17) \qquad \left(\frac{x_e}{x_{\tilde{e}}}, \frac{x_e}{x_v}, \frac{x_v}{x_{\tilde{v}}}\right)_{e,\tilde{e}\in E, v,\tilde{v}\in V} = \left(\frac{y_e}{y_{\tilde{e}}}, \frac{y_e}{y_v}, \frac{y_v}{y_{\tilde{v}}}\right)_{e,\tilde{e}\in E, v,\tilde{v}\in V}. \qquad \square$$

Part (b) of the lemma allows us to omit the reference edge $e_0$ in the normalizing constant:

$$z_{v_0} = z_{v_0,e_0}. \tag{3.18}$$

Hence, (3.12) can be rewritten as

$$\frac{d\mathbb{Q}_{v_0,e_0}}{d\rho_{e_0}}(x) = \frac{1}{z_{v_0}} \Phi_{v_0,\mathbf{a}}(x), \qquad x \in \Omega_{e_0}. \tag{3.19}$$

## 4. Interpolation between two environments.

For two vertices $v_0, v_1 \in V$, we consider an interpolation of the random environments $\mathbb{Q}_{v_0,e_0}$ and $\mathbb{Q}_{v_1,e_0}$ associated with the edge-reinforced random walk on $G$ with initial weights $\mathbf{a}$ and starting points $v_0$ and $v_1$, respectively:

DEFINITION 4.1 (Interpolated measure). For $v_0, v_1 \in V$ and $e_0 \in E$, define $\mathbb{P}_{v_0,v_1,e_0}$ to be the probability measure on $\Omega_{e_0}$ with the following Radon–Nikodym derivative with respect to $\rho_{e_0}$:

$$\frac{d\mathbb{P}_{v_0,v_1,e_0}}{d\rho_{e_0}} := \frac{1}{Z_{v_0,v_1,e_0}} \left( \frac{d\mathbb{Q}_{v_0,e_0}}{d\rho_{e_0}} \frac{d\mathbb{Q}_{v_1,e_0}}{d\rho_{e_0}} \right)^{1/2} \tag{4.1}$$

with the normalizing constant

$$Z_{v_0,v_1,e_0} := \int_{\Omega_{e_0}} \left( \frac{d\mathbb{Q}_{v_0,e_0}}{d\rho_{e_0}} \frac{d\mathbb{Q}_{v_1,e_0}}{d\rho_{e_0}} \right)^{1/2} d\rho_{e_0}. \tag{4.2}$$

REMARK 4.2. Since $d\mathbb{Q}_{v_0,e_0}/d\rho_{e_0}$ and $d\mathbb{Q}_{v_1,e_0}/d\rho_{e_0}$ are strictly positive probability densities with respect to $\rho_{e_0}$, the normalizing constant $Z_{v_0,v_1,e_0}$ is indeed strictly positive and, by the Cauchy–Schwarz inequality, it is finite. Hence, $\mathbb{P}_{v_0,v_1,e_0}$ is well defined.

The measure $\mathbb{P}_{v_0,v_1,e_0}$ inherits the following invariance under change of normalization from $\mathbb{Q}_{v_0,e_0}$:

LEMMA 4.3.

(a) *The normalizing constant $Z_{v_0,v_1,e_0}$ does not depend on the choice of the reference edge:*

$$Z_{v_0,v_1,e_0} = Z_{v_0,v_1,e_1} \qquad \text{for all } e_0, e_1 \in E. \tag{4.3}$$

(b) *The law of $g_{e_0,e_1}$, that is, the distribution of $(x_e/x_{e_1})_{e \in E}$, under $\mathbb{P}_{v_0,v_1,e_0}$ equals $\mathbb{P}_{v_0,v_1,e_1}$.*
(c) *Consequently, the ratios $(x_e/x_v)_{e \in E, v \in V}$ have the same joint laws under $\mathbb{P}_{v_0,v_1,e_0}$ and $\mathbb{P}_{v_0,v_1,e_1}$.*

PROOF. Let $e_0, e_1 \in E$. We claim that for any measurable function $f:\Omega_{e_1} \to [0, \infty]$ the following identity holds:

$$
\begin{aligned}
(4.4) \quad & \int_{\Omega_{e_0}} f(g_{e_0,e_1}(x))\left(\frac{d\mathbb{Q}_{v_0,e_0}}{d\rho_{e_0}}(x)\frac{d\mathbb{Q}_{v_1,e_0}}{d\rho_{e_0}}(x)\right)^{1/2}\rho_{e_0}(dx) \\
& = \int_{\Omega_{e_1}} f(x)\left(\frac{d\mathbb{Q}_{v_0,e_1}}{d\rho_{e_1}}(x)\frac{d\mathbb{Q}_{v_1,e_1}}{d\rho_{e_1}}(x)\right)^{1/2}\rho_{e_1}(dx).
\end{aligned}
$$

To see this, we first insert the description (3.19) of the random environments $\mathbb{Q}_{v_0,e_0}$ and $\mathbb{Q}_{v_1,e_0}$. Then, we apply the scaling invariance (3.9) of the densities, and finally, we use Lemma 3.1. This yields the following for the left-hand side of (4.4):

$$
\begin{aligned}
& \int_{\Omega_{e_0}} f(g_{e_0,e_1}(x))\left(\frac{1}{z_{v_0}}\Phi_{v_0,\mathbf{a}}(x)\frac{1}{z_{v_1}}\Phi_{v_1,\mathbf{a}}(x)\right)^{1/2}\rho_{e_0}(dx) \\
(4.5) \quad & = \int_{\Omega_{e_0}} f(g_{e_0,e_1}(x))\left(\frac{1}{z_{v_0}}\Phi_{v_0,\mathbf{a}}(g_{e_0,e_1}(x))\frac{1}{z_{v_1}}\Phi_{v_1,\mathbf{a}}(g_{e_0,e_1}(x))\right)^{1/2}\rho_{e_0}(dx) \\
& = \int_{\Omega_{e_1}} f(x)\left(\frac{1}{z_{v_0}}\Phi_{v_0,\mathbf{a}}(x)\frac{1}{z_{v_1}}\Phi_{v_1,\mathbf{a}}(x)\right)^{1/2}\rho_{e_1}(dx),
\end{aligned}
$$

which equals the right-hand side of (4.4).

The claim (4.3) in part (a) is a special case of (4.4) with $f \equiv 1$. Part (b) is an immediate consequence of part (a) and (4.4). Part (c) follows from part (b) in the same way as Lemma 3.5(d) follows from Lemma 3.5(c). □

Part (a) of the lemma allows us to drop the reference edge $e_0$ in the normalizing constant:

$$
(4.6) \qquad Z_{v_0,v_1} = Z_{v_0,v_1,e_0}.
$$

One reason why we introduce $\mathbb{P}_{v_0,v_1,e_0}$ is the following symmetry property, which does not hold for the random environment distribution $\mathbb{Q}_{v_0,e_0}$:

LEMMA 4.4 (Symmetry property of $\mathbb{P}_{v_0,v_1,e_0}$). *Let $v_0, v_1 \in V$ be such that Assumption* 2.3(c) *holds. Then, for all $e_0 \in E$, the ratios $x_{v_0}/x_{v_1}$ and $x_{v_1}/x_{v_0}$ have the same law under $\mathbb{P}_{v_0,v_1,e_0}$.*

PROOF. Let $v_0, v_1$ and $f$ be as in Assumption 2.3(c), and let $e_0 \in E$. We claim that the two random vectors

$$
(4.7) \quad \left(\frac{x_e}{x_{v_0}}\right)_{e\in E} \text{ and } \left(\frac{x_{f(e)}}{x_{v_1}}\right)_{e\in E} \text{ have the same law under } \mathbb{P}_{v_0,v_1,e_0}.
$$

To see this, let $A \subseteq \Omega$ be a measurable set. Using the definition of $\mathbb{P}_{v_0,v_1,e_0}$ and the description (3.19) of the random environment in the first step, we get

$$
\begin{aligned}
&\mathbb{P}_{v_0,v_1,e_0}\left[\left(\frac{x_e}{x_{v_0}}\right)_{e \in E} \in A\right] \\
&= \frac{1}{Z_{v_0,v_1}} \int_{\Omega_{e_0}} 1_A\left(\left(\frac{x_e}{x_{v_0}}\right)_{e \in E}\right)\left(\frac{1}{z_{v_0}}\Phi_{v_0,\mathbf{a}}(x)\frac{1}{z_{v_1}}\Phi_{v_1,\mathbf{a}}(x)\right)^{1/2} \rho_{e_0}(dx) \\
&= \frac{1}{Z_{v_0,v_1}} \int_{\Omega_{f(e_0)}} 1_A\left(\left(\frac{x_{f(e)}}{x_{f(v_0)}}\right)_{e \in E}\right) \\
&\quad \times \left(\frac{1}{z_{v_0}}\Phi_{v_0,\mathbf{a}}((x_{f(e)})_{e \in E})\frac{1}{z_{v_1}}\Phi_{v_1,\mathbf{a}}((x_{f(e)})_{e \in E})\right)^{1/2} \rho_{f(e_0)}(dx),
\end{aligned}
$$
(4.8)

where in the last identity we used that $\rho_{e_0}$ is obtained from $\rho_{f(e_0)}$ by permuting the variables by $x \mapsto (x_{f(e)})_{e \in E}$. Since $f$ is an automorphism of the weighted graph $(G, \mathbf{a})$ with $f(v_1) = v_0$, it follows from the description (3.7) of the densities that

$$
(4.9) \qquad \Phi_{v_0,\mathbf{a}}((x_{f(e)})_{e \in E}) = \Phi_{f^{-1}(v_0),\mathbf{a}}(x) = \Phi_{v_1,\mathbf{a}}(x).
$$

Using $f(v_0) = v_1$, one obtains in the same way $\Phi_{v_1,\mathbf{a}}((x_{f(e)})_{e \in E}) = \Phi_{v_0,\mathbf{a}}(x)$. Hence, we conclude from (4.8)

$$
\begin{aligned}
&\mathbb{P}_{v_0,v_1,e_0}\left[\left(\frac{x_e}{x_{v_0}}\right)_{e \in E} \in A\right] \\
&= \frac{1}{Z_{v_0,v_1}} \int_{\Omega_{f(e_0)}} 1_A\left(\left(\frac{x_{f(e)}}{x_{v_1}}\right)_{e \in E}\right) \\
&\quad \times \left(\frac{1}{z_{v_0}}\Phi_{v_1,\mathbf{a}}(x)\frac{1}{z_{v_1}}\Phi_{v_0,\mathbf{a}}(x)\right)^{1/2} \rho_{f(e_0)}(dx) \\
&= \mathbb{P}_{v_0,v_1,f(e_0)}\left[\left(\frac{x_{f(e)}}{x_{v_1}}\right)_{e \in E} \in A\right].
\end{aligned}
$$
(4.10)

Lemma 4.3(c) implies that in the last probability one can replace the reference edge $f(e_0)$ by $e_0$. Hence, the claim (4.7) follows.

In particular, we conclude from (4.7) that

$$
(4.11) \qquad \sum_{e \ni v_1} \frac{x_e}{x_{v_0}} = \frac{x_{v_1}}{x_{v_0}} \quad \text{and} \quad \sum_{e \ni v_1} \frac{x_{f(e)}}{x_{v_1}} = \frac{x_{v_0}}{x_{v_1}}
$$

have the same law under $\mathbb{P}_{v_0,v_1,e_0}$. □

Next, we present the connection between the interpolated measure and the random environment distribution. The normalizing constant $Z_{v_0,v_1}$ of the interpolated measure turns out to be equal to the 1/4-moment estimated in the main lemma 2.5; see the left-hand side of (2.19).

LEMMA 4.5 (Relation between $\mathbb{P}_{v_0,v_1,e_0}$ and $\mathbb{Q}_{v_0,e_0}$). *Let $v_0, v_1 \in V$ be such that Assumption 2.3(c) holds. Then, for all $e_0 \in E$, the measures $\mathbb{P}_{v_0,v_1,e_0}$ and $\mathbb{Q}_{v_0,e_0}$ are mutually absolutely continuous with the Radon–Nikodym derivative*

$$(4.12) \qquad \frac{d\mathbb{P}_{v_0,v_1,e_0}}{d\mathbb{Q}_{v_0,e_0}} = \frac{1}{Z_{v_0,v_1}} \left(\frac{x_{v_1}}{x_{v_0}}\right)^{1/4}.$$

*As a consequence, one has*

$$(4.13) \qquad Z_{v_0,v_1} = E_{\mathbb{Q}_{v_0,e_0}}\left[\left(\frac{x_{v_1}}{x_{v_0}}\right)^{1/4}\right].$$

PROOF. Let $v_0, v_1$ and $f$ be as in Assumption 2.3(c), and let $e_0 \in E$. Using the definition of $\mathbb{P}_{v_0,v_1,e_0}$ and formulas (3.19) and (3.7) for the densities $d\mathbb{Q}_{v_1,e_0}/d\rho_{e_0}$ and $d\mathbb{Q}_{v_0,e_0}/d\rho_{e_0}$, we get

$$(4.14) \qquad \begin{aligned} \frac{d\mathbb{P}_{v_0,v_1,e_0}}{d\mathbb{Q}_{v_0,e_0}} &= \frac{1}{Z_{v_0,v_1}} \left(\frac{d\mathbb{Q}_{v_1,e_0}}{d\rho_{e_0}}\right)^{1/2} \left(\frac{d\mathbb{Q}_{v_0,e_0}}{d\rho_{e_0}}\right)^{-1/2} \\ &= \frac{1}{Z_{v_0,v_1}} \left(\frac{z_{v_0}}{z_{v_1}}\right)^{1/2} \left(\frac{x_{v_1}}{x_{v_0}}\right)^{1/4}. \end{aligned}$$

We claim that $z_{v_1} = z_{v_0}$. Indeed, since $f$ is an automorphism of the weighted graph $(G, \mathbf{a})$ with $f(v_0) = v_1$, we can use (4.9) to obtain

$$(4.15) \qquad \begin{aligned} z_{v_1} &= \int_{\Omega_{e_0}} \Phi_{v_1,\mathbf{a}}(x)\rho_{e_0}(dx) = \int_{\Omega_{e_0}} \Phi_{v_0,\mathbf{a}}((x_{f(e)})_{e \in E})\rho_{e_0}(dx) \\ &= \int_{\Omega_{f^{-1}(e_0)}} \Phi_{v_0,\mathbf{a}}(x)\rho_{f^{-1}(e_0)}(dx) = z_{v_0}; \end{aligned}$$

recall that these normalizing constants do not depend on the choice of the reference edge $e_0$ by Lemma 3.5(b). Thus, (4.12) follows from (4.14).

Integrating (4.12) with respect to the probability measure $\mathbb{Q}_{v_0,e_0}$ yields the second claim. □

REMARK 4.6. The fact $Z_{v_0,v_1} < \infty$, which was noted in Remark 4.2, and the identity (4.13) together imply $E_{\mathbb{Q}_{v_0,e_0}}[(x_{v_1}/x_{v_0})^{1/4}] < \infty$.

**5. The key estimate in the proof of Lemma 2.5.** From now on until the end of Section 7, we assume that Assumption 2.3 holds.

In order to prove the main Lemma 2.5, it suffices by (4.13) to get good estimates for the normalizing constant $Z_{v_0,v_1}$. In this section, we derive such an upper bound for $Z_{v_0,v_1}$.

We introduce the following function of the random environment:

$$(5.1) \qquad H_{v_0,v_1} : \Omega \to \mathbb{R}, \qquad H_{v_0,v_1}(x) := \frac{1}{4} \log \frac{x_{v_1}}{x_{v_0}}.$$

Lemma 4.5 allows us to write

$$(5.2) \qquad \mathbb{P}_{v_0,v_1,e_0}(dx) = \frac{e^{H_{v_0,v_1}(x)}}{Z_{v_0,v_1}} \mathbb{Q}_{v_0,e_0}(dx).$$

The following general variational principle is a cornerstone of this paper.

Let $\mathcal{P}$ denote the set of all probability measures $\Pi$ on $\Omega_{e_0}$ such that $\Pi$ and $\mathbb{P}_{v_0,v_1,e_0}$ are mutually absolutely continuous and $E_\Pi[H_{v_0,v_1}]$ and $E_\Pi[\log(d\Pi/d\mathbb{P}_{v_0,v_1,e_0})]$ are finite.

LEMMA 5.1 (Variational principle). *The functional*

$$(5.3) \qquad F : \mathcal{P} \to \mathbb{R}, \qquad F(\Pi) := E_\Pi[H_{v_0,v_1}] + E_\Pi\left[\log \frac{d\Pi}{d\mathbb{P}_{v_0,v_1,e_0}}\right]$$

*is minimized for* $\Pi = \mathbb{Q}_{v_0,e_0} \in \mathcal{P}$ *with the value*

$$(5.4) \qquad F(\mathbb{Q}_{v_0,e_0}) = \log Z_{v_0,v_1}.$$

PROOF. We prove first that $\mathbb{Q}_{v_0,e_0} \in \mathcal{P}$. By Lemma 4.5, $\mathbb{Q}_{v_0,e_0}$ and $\mathbb{P}_{v_0,v_1,e_0}$ are mutually absolutely continuous. To see that the expectation $E_{\mathbb{Q}_{v_0,e_0}}[H_{v_0,v_1}]$ is finite, we use Theorem 2.4 in [8]. In the present context, this theorem states that for all edges $e, e' \in E$, with respect to $\mathbb{Q}_{v_0,e_0}$, the random variables $x_e/x_{e'}$ have tails bounded by a power law. As a consequence, the law of $H_{v_0,v_1}(x) = 4^{-1} \log(x_{v_1}/x_{v_0})$ has exponentially bounded tails, that is,

$$(5.5) \qquad \limsup_{M \to \infty} \frac{1}{M} \log \mathbb{Q}_{v_0,e_0}\left[\left|\log \frac{x_{v_1}}{x_{v_0}}\right| \geq M\right] < 0.$$

This implies $E_{\mathbb{Q}_{v_0,e_0}}[|H_{v_0,v_1}|] < \infty$.

Next, we show that $E_{\mathbb{Q}_{v_0,e_0}}[\log(d\mathbb{Q}_{v_0,e_0}/d\mathbb{P}_{v_0,v_1,e_0})]$ is finite. The representation (5.2) of $\mathbb{P}_{v_0,v_1,e_0}$ implies

$$(5.6) \qquad \log \frac{d\mathbb{Q}_{v_0,e_0}}{d\mathbb{P}_{v_0,v_1,e_0}} = \log Z_{v_0,v_1} - H_{v_0,v_1}.$$

From Remark 4.2 we know that $\log Z_{v_0,v_1}$ is finite. Using (5.6) and $E_{\mathbb{Q}_{v_0,e_0}}[|H_{v_0,v_1}|] < \infty$, this implies that $\log(d\mathbb{Q}_{v_0,e_0}/d\mathbb{P}_{v_0,v_1,e_0})$ has a finite expectation with respect to $\mathbb{Q}_{v_0,e_0}$. This completes the proof of $\mathbb{Q}_{v_0,e_0} \in \mathcal{P}$.

Taking the expectation of (5.6) with respect to $\mathbb{Q}_{v_0,e_0}$ yields the claim $F(\mathbb{Q}_{v_0,e_0}) = \log Z_{v_0,v_1}$ as follows:

$$(5.7) \quad \log Z_{v_0,v_1} = E_{\mathbb{Q}_{v_0,e_0}}[H_{v_0,v_1}] + E_{\mathbb{Q}_{v_0,e_0}}\left[\log \frac{d\mathbb{Q}_{v_0,e_0}}{d\mathbb{P}_{v_0,v_1,e_0}}\right] = F(\mathbb{Q}_{v_0,e_0}).$$

We now prove that $\log Z_{v_0,v_1}$ is indeed the minimal value of the functional $F$. Let $\Pi \in \mathcal{P}$. We rewrite (5.6) in the following form:

$$(5.8) \quad \begin{aligned} \log Z_{v_0,v_1} &= H_{v_0,v_1} + \log \frac{d\mathbb{Q}_{v_0,e_0}}{d\mathbb{P}_{v_0,v_1,e_0}} \\ &= H_{v_0,v_1} + \log \frac{d\Pi}{d\mathbb{P}_{v_0,v_1,e_0}} - \log \frac{d\Pi}{d\mathbb{Q}_{v_0,e_0}}. \end{aligned}$$

Taking expectations with respect to $\Pi$, this implies

$$(5.9) \quad \log Z_{v_0,v_1} = F(\Pi) - E_\Pi\left[\log \frac{d\Pi}{d\mathbb{Q}_{v_0,e_0}}\right];$$

recall the definition (5.3) of $F$. Now, the last relative entropy is nonnegative:

$$(5.10) \quad E_\Pi\left[\log \frac{d\Pi}{d\mathbb{Q}_{v_0,e_0}}\right] \geq 0.$$

This implies the claim $F(\Pi) \geq \log Z_{v_0,v_1}$. $\square$

Before going into the details, let us intuitively and roughly explain how to get a good upper bound for $\log Z_{v_0,v_1}$. To apply the variational principle, we need to construct an appropriate $\Pi = \Pi_{v_0,v_1,e_0}^{(\gamma\varphi)} \in \mathcal{P}$ depending on a parameter $\gamma \in \mathbb{R}$. On the one hand, $\Pi$ should be "close" to $\mathbb{Q}_{v_0,e_0}$ to get the relative entropy in (5.10) small. On the other hand, we need to control $F(\Pi)$. This requires us to control both, $E_\Pi[H_{v_0,v_1}]$ and the relative entropy $E_\Pi[\log(d\Pi/d\mathbb{P}_{v_0,v_1,e_0})]$. We choose $\Pi$ to be a "deformation" of $\mathbb{P}_{v_0,v_1,e_0}$. By this we mean that $\Pi$ is the law of a suitable random variable $\Xi_{e_0}^{(\gamma\varphi)} : \Omega_{e_0} \to \Omega_{e_0}$ with respect to $\mathbb{P}_{v_0,v_1,e_0}$. We call $\Xi_{e_0}^{(\gamma\varphi)}$ the "deformation map." It depends on a deformation parameter $\gamma \in \mathbb{R}$.

The intuition behind our choice of the deformation map $\Xi_{e_0}^{(\gamma\varphi)}$ is as follows: In the case of large boxes $G_r^{(i)}$, for $v_0 = \mathbf{0}$ and $v_1 \in L_r^{(i)}$ far away from the origin, roughly speaking, one expects that under $\mathbb{Q}_{v_0,e_0}$, the weights of the edges close to $v_1$ are typically much smaller than the weights of the edges close to $v_0$. On the other hand, under the interpolated measure $\mathbb{P}_{v_0,v_1,e_0}$, one might expect intuitively that the edge weights close to $v_0$ and close to $v_1$ are more or less of the same order of magnitude. This motivates us to define $\Xi_{e_0}^{(\gamma\varphi)}$ below in such a way that it leaves the weights of the edges incident to $v_0$ unchanged [see (5.14) below] and multiplies the weights of the edges incident to $v_1$ by $e^\gamma$ [see (5.15) below]. For negative $\gamma$, this yields the desired behavior.

Formally, the deformation map is defined as follows:

DEFINITION 5.2. Let $v_0, v_1, e_0$ and $\varphi$ be as in Assumption 2.3. For $\gamma \in \mathbb{R}$, we define the map

$$\Xi^{(\gamma\varphi)} : \Omega \to \Omega, \qquad (x_e)_{e \in E} \mapsto (e^{\gamma\varphi(e)} x_e)_{e \in E}. \tag{5.11}$$

As a consequence of Assumptions 2.3(b) and (d), $(\Xi^{(\gamma\varphi)}(x))_{e_0} = x_{e_0}$ and, hence, $\Xi^{(\gamma\varphi)}$ maps $\Omega_{e_0}$ to $\Omega_{e_0}$. We denote the restriction of $\Xi^{(\gamma\varphi)}$ to $\Omega_{e_0}$ by

$$\Xi^{(\gamma\varphi)}_{e_0} : \Omega_{e_0} \to \Omega_{e_0}. \tag{5.12}$$

Let $\Pi^{(\gamma\varphi)}_{v_0,v_1,e_0}$ denote the distribution of $\Xi^{(\gamma\varphi)}_{e_0}$ under $\mathbb{P}_{v_0,v_1,e_0}$.

Let us state two important properties of the map $\Xi^{(\gamma\varphi)}_{e_0}$:

LEMMA 5.3. *For all $\gamma \in \mathbb{R}$ and $x \in \Omega_{e_0}$, one has*

$$(\Xi^{(\gamma\varphi)}_{e_0}(x))_{v_0} = x_{v_0} \quad and \quad (\Xi^{(\gamma\varphi)}_{e_0}(x))_{v_1} = e^\gamma x_{v_1}. \tag{5.13}$$

PROOF. Let $x \in \Omega_{e_0}$. It follows from the properties of $\varphi$ assumed in Assumptions 2.3(d) that

$$(\Xi^{(\gamma\varphi)}_{e_0}(x))_e = x_e \qquad \text{for all } e \ni v_0 \quad \text{and} \tag{5.14}$$

$$(\Xi^{(\gamma\varphi)}_{e_0}(x))_e = e^\gamma x_e \qquad \text{for all } e \ni v_1. \tag{5.15}$$

Summing (5.14) over all edges $e$ incident to $v_0$ and (5.15) over all edges $e$ incident to $v_1$ yields the claims in (5.13). $\square$

Applying the variational principle (Lemma 5.1) with the measures $\Pi^{(\gamma\varphi)}_{v_0,v_1,e_0}$, we obtain the following upper bound:

LEMMA 5.4 (Key estimate). *For all $\gamma \in \mathbb{R}$, one has*

$$\log Z_{v_0,v_1} \le E_{\Pi^{(\gamma\varphi)}_{v_0,v_1,e_0}}[H_{v_0,v_1}] + E_{\Pi^{(\gamma\varphi)}_{v_0,v_1,e_0}}\left[\log \frac{d\Pi^{(\gamma\varphi)}_{v_0,v_1,e_0}}{d\mathbb{P}_{v_0,v_1,e_0}}\right]. \tag{5.16}$$

To prove this key estimate, we need to verify that $\Pi^{(\gamma\varphi)}_{v_0,v_1,e_0} \in \mathcal{P}$ for all $\gamma \in \mathbb{R}$. This requires some preparations.

We start by introducing additional notation. For $\gamma \in \mathbb{R}$ and $x \in \Omega$, we set

$$(5.17) \quad x_v^{(\gamma\varphi)} := \sum_{\substack{e \in E: \\ e \ni v}} e^{\gamma\varphi(e)} x_e \quad \text{for every vertex } v \in V,$$

$$(5.18) \quad Y_T^{(\gamma\varphi)}(x) := \prod_{e \in T} (e^{\gamma\varphi(e)} x_e) \quad \text{for every spanning tree } T \in \mathcal{T}.$$

In the special case $\gamma = 0$, we have $x_v^{(0\varphi)} = x_v$ and we abbreviate

$$(5.19) \quad Y_T(x) := Y_T^{(0\varphi)}(x) = \prod_{e \in T} x_e.$$

LEMMA 5.5. *For all $\gamma \in \mathbb{R}$, the measures $\Pi_{v_0,v_1,e_0}^{(\gamma\varphi)}$ and $\mathbb{P}_{v_0,v_1,e_0}$ are mutually absolutely continuous. The Radon–Nikodym derivative satisfies the following identity for $\rho_{e_0}$-almost all $x = (x_e)_{e \in E} \in \Omega_{e_0}$:*

$$(5.20) \quad \begin{aligned} &\log \frac{d\Pi_{v_0,v_1,e_0}^{(\gamma\varphi)}}{d\mathbb{P}_{v_0,v_1,e_0}}(\Xi_{e_0}^{(\gamma\varphi)}(x)) \\ &= \gamma\left(\frac{a_{v_1}}{2} + \frac{1}{4}\right) - \sum_{e \in E} \gamma\varphi(e) a_e \\ &\quad + \sum_{v \in V \setminus \{v_0,v_1\}} \frac{a_v + 1}{2} \log \frac{x_v^{(\gamma\varphi)}}{x_v} - \frac{1}{2} \log \frac{\sum_{T \in \mathcal{T}} Y_T^{(\gamma\varphi)}(x)}{\sum_{T \in \mathcal{T}} Y_T(x)}. \end{aligned}$$

In the remainder of this article, we work with this particular realization of the Radon–Nikodym derivative $d\Pi_{v_0,v_1,e_0}^{(\gamma\varphi)}/d\mathbb{P}_{v_0,v_1,e_0}$.

PROOF OF LEMMA 5.5. By its definition, $\mathbb{P}_{v_0,v_1,e_0}$ is absolutely continuous with respect to $\rho_{e_0}$. Consequently, since $\Pi_{v_0,v_1,e_0}^{(\gamma\varphi)}$ is the law of $\Xi_{e_0}^{(\gamma\varphi)}$ under $\mathbb{P}_{v_0,v_1,e_0}$ due to Definition 5.2, one has for any positive measurable function $h: \Omega_{e_0} \to \mathbb{R}$,

$$(5.21) \quad \begin{aligned} E_{\Pi_{v_0,v_1,e_0}^{(\gamma\varphi)}}[h] &= E_{\mathbb{P}_{v_0,v_1,e_0}}[h \circ \Xi_{e_0}^{(\gamma\varphi)}] \\ &= \int_{\Omega_{e_0}} h(\Xi_{e_0}^{(\gamma\varphi)}(x)) \frac{d\mathbb{P}_{v_0,v_1,e_0}}{d\rho_{e_0}}(x) \rho_{e_0}(dx). \end{aligned}$$

The map $\Xi_{e_0}^{(\gamma\varphi)}$ is invertible. Due to (5.14), $\Xi_{e_0}^{(\gamma\varphi)}$ leaves the weight $x_{e_0}$ of the reference edge $e_0 \ni v_0$ unchanged and it multiplies all other $x_e$ with some $e$-dependent constant. Thus, $\Xi_{e_0}^{(\gamma\varphi)}$ leaves the reference measure $\rho_{e_0}$ invariant: $\rho_{e_0}[[\Xi_{e_0}^{(\gamma\varphi)}]^{-1} A] = \rho_{e_0}[A]$ for all measurable sets $A \subseteq \Omega_{e_0}$. Consequently, applying the transformation formula on the right-hand side of (5.21), we obtain

$$(5.22) \qquad E_{\Pi^{(\gamma\varphi)}_{v_0,v_1,e_0}}[h] = \int_{\Omega_{e_0}} h(x) \frac{d\mathbb{P}_{v_0,v_1,e_0}}{d\rho_{e_0}}([\Xi^{(\gamma\varphi)}_{e_0}]^{-1}(x))\rho_{e_0}(dx).$$

This implies that $\Pi^{(\gamma\varphi)}_{v_0,v_1,e_0}$ is absolutely continuous with respect to $\rho_{e_0}$ with the Radon–Nikodym derivative

$$(5.23) \qquad \frac{d\Pi^{(\gamma\varphi)}_{v_0,v_1,e_0}}{d\rho_{e_0}} = \frac{d\mathbb{P}_{v_0,v_1,e_0}}{d\rho_{e_0}} \circ [\Xi^{(\gamma\varphi)}_{e_0}]^{-1}.$$

Using (3.12) and Lemma 4.5, it follows that both $d\mathbb{P}_{v_0,v_1,e_0}/d\rho_{e_0}$ and $d\Pi^{(\gamma\varphi)}_{v_0,v_1,e_0}/d\rho_{e_0}$ are strictly positive on $\Omega_{e_0}$. Hence, the measures $\Pi^{(\gamma\varphi)}_{v_0,v_1,e_0}$ and $\mathbb{P}_{v_0,v_1,e_0}$ are mutually absolutely continuous, and for $\rho_{e_0}$-almost all $x \in \Omega_{e_0}$, the Radon–Nikodym derivative satisfies

$$(5.24) \quad \begin{aligned} &\log \frac{d\Pi^{(\gamma\varphi)}_{v_0,v_1,e_0}}{d\mathbb{P}_{v_0,v_1,e_0}}(\Xi^{(\gamma\varphi)}_{e_0}(x)) \\ &= \log\left(\frac{d\Pi^{(\gamma\varphi)}_{v_0,v_1,e_0}}{d\rho_{e_0}}(\Xi^{(\gamma\varphi)}_{e_0}(x)) \cdot \left[\frac{d\mathbb{P}_{v_0,v_1,e_0}}{d\rho_{e_0}}(\Xi^{(\gamma\varphi)}_{e_0}(x))\right]^{-1}\right) \\ &= \log\left(\frac{d\mathbb{P}_{v_0,v_1,e_0}}{d\rho_{e_0}}(x) \cdot \left[\frac{d\mathbb{P}_{v_0,v_1,e_0}}{d\rho_{e_0}}(\Xi^{(\gamma\varphi)}_{e_0}(x))\right]^{-1}\right). \end{aligned}$$

To get an explicit form for the last expression, we use Lemma 4.5 and (5.13) to obtain

$$(5.25) \qquad \frac{d\mathbb{P}_{v_0,v_1,e_0}}{d\rho_{e_0}}(\Xi^{(\gamma\varphi)}_{e_0}(x)) = \frac{1}{Z_{v_0,v_1}}\left(\frac{e^\gamma x_{v_1}}{x_{v_0}}\right)^{1/4} \frac{d\mathbb{Q}_{v_0,e_0}}{d\rho_{e_0}}(\Xi^{(\gamma\varphi)}_{e_0}(x)).$$

Next, we insert formulas (3.19) and (3.7) for $d\mathbb{Q}_{v_0,e_0}/d\rho_{e_0}$ and simplify the resulting expression using the abbreviations (5.17) and (5.18):

$$(5.26) \quad \begin{aligned} &\frac{d\mathbb{P}_{v_0,v_1,e_0}}{d\rho_{e_0}}(\Xi^{(\gamma\varphi)}_{e_0}(x)) \\ &= \frac{1}{Z_{v_0,v_1}z_{v_0}} \\ &\quad \times \frac{\prod_{e\in E}(e^{\gamma\varphi(e)a_e}x_e^{a_e})\sqrt{\sum_{T\in\mathcal{T}}\prod_{e\in T}(e^{\gamma\varphi(e)}x_e)}}{x_{v_0}^{a_{v_0}/2+1/4}(e^\gamma x_{v_1})^{a_{v_1}/2+1/4}\prod_{v\in V\setminus\{v_0,v_1\}}[\sum_{e\ni v}e^{\gamma\varphi(e)}x_e]^{(a_v+1)/2}} \\ &= \frac{1}{Z_{v_0,v_1}z_{v_0}} \end{aligned}$$

$$\times \frac{\prod_{e \in E} e^{\gamma \varphi(e) a_e}}{e^{\gamma(a_{v_1}/2+1/4)}} \frac{\prod_{e \in E} x_e^{a_e} \sqrt{\sum_{T \in \mathcal{T}} Y_T^{(\gamma\varphi)}(x)}}{x_{v_0}^{a_{v_0}/2+1/4} x_{v_1}^{a_{v_1}/2+1/4} \prod_{v \in V \setminus \{v_0, v_1\}} (x_v^{(\gamma\varphi)})^{(a_v+1)/2}}.$$

For $\gamma = 0$, $\Xi_{e_0}^{(\gamma\varphi)} \colon \Omega_{e_0} \to \Omega_{e_0}$ is the identity map by Definition 5.2: $\Xi_{e_0}^{(\gamma\varphi)}(x) = x$ for all $x$. Hence, we obtain

$$(5.27) \quad \begin{aligned} &\frac{d\mathbb{P}_{v_0,v_1,e_0}}{d\rho_{e_0}}(x) \cdot \left[\frac{d\mathbb{P}_{v_0,v_1,e_0}}{d\rho_{e_0}}(\Xi_{e_0}^{(\gamma\varphi)}(x))\right]^{-1} \\ &= \frac{e^{\gamma(a_{v_1}/2+1/4)}}{\prod_{e \in E} e^{\gamma\varphi(e)a_e}} \prod_{v \in V \setminus \{v_0,v_1\}} \left(\frac{x_v^{(\gamma\varphi)}}{x_v}\right)^{(a_v+1)/2} \sqrt{\frac{\sum_{T \in \mathcal{T}} Y_T(x)}{\sum_{T \in \mathcal{T}} Y_T^{(\gamma\varphi)}(x)}}. \end{aligned}$$

Taking logarithms in the last formula, the claim follows from (5.24). $\square$

LEMMA 5.6. *For all $\gamma \in \mathbb{R}$, one has*

$$(5.28) \qquad 0 \le E_{\Pi_{v_0,v_1,e_0}^{(\gamma\varphi)}}\left[\log \frac{d\Pi_{v_0,v_1,e_0}^{(\gamma\varphi)}}{d\mathbb{P}_{v_0,v_1,e_0}}\right] < \infty.$$

PROOF. By Assumption 2.3, $\varphi$ takes only values in $[0,1]$, so that $e^{-|\gamma|} \le e^{\gamma\varphi(e)} \le e^{|\gamma|}$ for all edges $e$ and all $\gamma \in \mathbb{R}$. This bound implies $e^{-|\gamma|}x_v \le x_v^{(\gamma\varphi)} \le e^{|\gamma|}x_v$ and $e^{-|\gamma|\cdot(|V|-1)}Y_T(x) \le Y_T^{(\gamma\varphi)}(x) \le e^{|\gamma|\cdot(|V|-1)}Y_T(x)$ for all $v \in V$ and all $T \in \mathcal{T}$; here we use that every spanning tree $T \in \mathcal{T}$ consists of $|V|-1$ edges. Thus, using (5.20) and the fact that $\Xi_{e_0}^{(\gamma\varphi)}$ is a bijection on $\Omega_{e_0}$, we get

$$(5.29) \quad \sup_{x \in \Omega_{e_0}} \left|\log \frac{d\Pi_{v_0,v_1,e_0}^{(\gamma\varphi)}}{d\mathbb{P}_{v_0,v_1,e_0}}(x)\right| = \sup_{x \in \Omega_{e_0}} \left|\log \frac{d\Pi_{v_0,v_1,e_0}^{(\gamma\varphi)}}{d\mathbb{P}_{v_0,v_1,e_0}}(\Xi_{e_0}^{(\gamma\varphi)}(x))\right| < \infty.$$

Hence, the expectation in (5.28) is finite. Note that this term is a relative entropy of $\Pi_{v_0,v_1,e_0}^{(\gamma\varphi)}$ with respect to $\mathbb{P}_{v_0,v_1,e_0}$; therefore, it is nonnegative. This completes the proof of the lemma. $\square$

In the proof of $\Pi_{v_0,v_1,e_0}^{(\gamma\varphi)} \in \mathcal{P}$, we need that the first expectation $E_{\Pi_{v_0,v_1,e_0}^{(\gamma\varphi)}}[H_{v_0,v_1}]$ arising in the upper bound of inequality (5.16) is finite. Later, for controlling the right-hand side of the key estimate (5.16), we need the stronger statement that this expectation is linear in $\gamma$:

LEMMA 5.7. *For all $\gamma \in \mathbb{R}$, one has*

$$(5.30) \qquad E_{\Pi_{v_0,v_1,e_0}^{(\gamma\varphi)}}[H_{v_0,v_1}] = \frac{\gamma}{4}.$$

ERRW ON A 2-DIMENSIONAL GRAPH 23PROOF. Recall the definition (5.1) of $H_{v_0,v_1}$. In the case $\gamma = 0$, one has $\Xi_{e_0}^{(\gamma\varphi)}(x) = x$ for all $x$ and, consequently, $\Pi_{v_0,v_1,e_0}^{(0\varphi)} = \mathbb{P}_{v_0,v_1,e_0}$. Hence, integrating (5.6) with respect to $\mathbb{P}_{v_0,v_1,e_0}$ yields

$$\frac{1}{4} E_{\mathbb{P}_{v_0,v_1,e_0}}\left[\log \frac{x_{v_1}}{x_{v_0}}\right] = E_{\mathbb{P}_{v_0,v_1,e_0}}[H_{v_0,v_1}]$$

$$= \log Z_{v_0,v_1} - E_{\mathbb{P}_{v_0,v_1,e_0}}\left[\log \frac{d\mathbb{Q}_{v_0,e_0}}{d\mathbb{P}_{v_0,v_1,e_0}}\right]$$

(5.31)

$$= \log Z_{v_0,v_1} + E_{\mathbb{P}_{v_0,v_1,e_0}}\left[\log \frac{d\mathbb{P}_{v_0,v_1,e_0}}{d\mathbb{Q}_{v_0,e_0}}\right]$$

$$\geq \log Z_{v_0,v_1} > -\infty.$$

Here we used that the relative entropy of $\mathbb{P}_{v_0,v_1,e_0}$ with respect to $\mathbb{Q}_{v_0,e_0}$ is non-negative. By Lemma 4.4, the random variables $\log(x_{v_1}/x_{v_0})$ and $\log(x_{v_0}/x_{v_1}) = -\log(x_{v_1}/x_{v_0})$ have the same law with respect to $\mathbb{P}_{v_0,v_1,e_0}$. Hence, all expectations in (5.31) are finite and, in particular, it follows that

(5.32) $$E_{\mathbb{P}_{v_0,v_1,e_0}}\left[\log \frac{x_{v_1}}{x_{v_0}}\right] = 0,$$

which proves the claim in the case $\gamma = 0$.

For $\gamma \in \mathbb{R}$, we use the definition of $\Pi_{v_0,v_1,e_0}^{(\gamma\varphi)}$ as the law of $\Xi_{e_0}^{(\gamma\varphi)}$ under $\mathbb{P}_{v_0,v_1,e_0}$ and the relations (5.13) to obtain

(5.33)
$$\frac{1}{4} E_{\Pi_{v_0,v_1,e_0}^{(\gamma\varphi)}}\left[\log \frac{x_{v_1}}{x_{v_0}}\right] = \frac{1}{4} E_{\mathbb{P}_{v_0,v_1,e_0}}\left[\log \frac{(\Xi_{e_0}^{(\gamma\varphi)}(x))_{v_1}}{(\Xi_{e_0}^{(\gamma\varphi)}(x))_{v_0}}\right]$$

$$= \frac{1}{4} E_{\mathbb{P}_{v_0,v_1,e_0}}\left[\log \frac{e^\gamma x_{v_1}}{x_{v_0}}\right]$$

$$= \frac{\gamma}{4} + \frac{1}{4} E_{\mathbb{P}_{v_0,v_1,e_0}}\left[\log \frac{x_{v_1}}{x_{v_0}}\right] = \frac{\gamma}{4}. \quad \square$$

Now we are ready to prove the key estimate:

PROOF OF LEMMA 5.4. Let $\gamma \in \mathbb{R}$. By the variational principle (Lemma 5.1), it only remains to verify that $\Pi_{v_0,v_1,e_0}^{(\gamma\varphi)} \in \mathcal{P}$.

By Lemma 5.5, the measures $\Pi_{v_0,v_1,e_0}^{(\gamma\varphi)}$ and $\mathbb{P}_{v_0,v_1,e_0}$ are mutually absolutely continuous. Lemma 5.7 shows that the expectation $E_{\Pi_{v_0,v_1,e_0}^{(\gamma\varphi)}}[H_{v_0,v_1}]$ is finite. The relative entropy $E_{\Pi_{v_0,v_1,e_0}^{(\gamma\varphi)}}[\log(d\Pi_{v_0,v_1,e_0}^{(\gamma\varphi)}/d\mathbb{P}_{v_0,v_1,e_0})]$ is finite by Lemma 5.6. This implies the claim $\Pi_{v_0,v_1,e_0}^{(\gamma\varphi)} \in \mathcal{P}$. $\quad \square$

**6. Estimating the relative entropy.** By Lemma 5.7, the first term in the key estimate (5.16) is linear in $\gamma$. The goal of this section is to prove the following quadratic bound in $\gamma$ for the second term in (5.16):

LEMMA 6.1 (Bound of the relative entropy). *Suppose that Assumption 2.3 holds. For all $\gamma \in \mathbb{R}$, one has*

$$(6.1) \qquad E_{\Pi_{v_0,v_1,e_0}^{(\gamma\varphi)}}\left[\log \frac{d\Pi_{v_0,v_1,e_0}^{(\gamma\varphi)}}{d\mathbb{P}_{v_0,v_1,e_0}}\right] \leq \frac{S_\varphi \gamma^2}{2}$$

*with $S_\varphi$ defined in (2.20).*

Before proving the lemma, let us introduce some notation. For $\gamma \in \mathbb{R}$, set

$$(6.2) \qquad g(\gamma) := E_{\Pi_{v_0,v_1,e_0}^{(\gamma\varphi)}}\left[\log \frac{d\Pi_{v_0,v_1,e_0}^{(\gamma\varphi)}}{d\mathbb{P}_{v_0,v_1,e_0}}\right].$$

By the definition of $\Pi_{v_0,v_1,e_0}^{(\gamma\varphi)}$ (Definition 5.2), we have $g(\gamma) = E_{\mathbb{P}_{v_0,v_1,e_0}}[f_\gamma]$ with

$$(6.3) \qquad f_\gamma(x) := \log \frac{d\Pi_{v_0,v_1,e_0}^{(\gamma\varphi)}}{d\mathbb{P}_{v_0,v_1,e_0}}(\Xi_{e_0}^{(\gamma\varphi)}(x))$$

for $x \in \Omega_{e_0}$. The idea for the bound of the relative entropy is to do a Taylor expansion of $g$ around $\gamma = 0$. First, we need to show that $g$ is twice continuously differentiable. This will be done by showing that $f_\gamma$ is twice continuously differentiable in $\gamma$ with derivatives which are bounded as functions of $x$.

Recall the abbreviations (5.17) and (5.18). It will be convenient to use the following notation: For $x \in \Omega_{e_0}$, $v \in V$, and $\gamma \in \mathbb{R}$, consider $e \mapsto \varphi(e)$ as a random variable on the space $\{e \in E : e \ni v\}$ equipped with the probability measure

$$(6.4) \qquad \mu_{x,v,\gamma} := \sum_{\substack{e \in E:\\ e \ni v}} \frac{e^{\gamma\varphi(e)} x_e}{x_v^{(\gamma\varphi)}} \delta_e.$$

For any spanning tree $T \in \mathcal{T}$, we abbreviate

$$(6.5) \qquad \Delta(T) := \sum_{e \in T} \varphi(e).$$

We consider the map $\Delta : T \mapsto \Delta(T)$ as a random variable on the space $\mathcal{T}$ of spanning trees equipped with the probability measure

$$(6.6) \qquad \nu_{x,\gamma} := \sum_{T \in \mathcal{T}} \frac{Y_T^{(\gamma\varphi)}(x)}{\sum_{T' \in \mathcal{T}} Y_{T'}^{(\gamma\varphi)}(x)} \delta_T.$$

LEMMA 6.2. *For each $x \in \Omega_{e_0}$, the function $\gamma \mapsto f_\gamma(x)$ is twice continuously differentiable with the following derivatives:*

$$\frac{\partial f_\gamma}{\partial \gamma}(x) = \frac{a_{v_1}}{2} + \frac{1}{4} - \sum_{e \in E} \varphi(e) a_e$$
(6.7)
$$+ \sum_{v \in V \setminus \{v_0, v_1\}} \frac{a_v + 1}{2} E_{\mu_{x,v,\gamma}}[\varphi] - \frac{1}{2} E_{\nu_{x,\gamma}}[\Delta],$$

(6.8) $$\frac{\partial^2 f_\gamma}{\partial^2 \gamma}(x) = \sum_{v \in V \setminus \{v_0, v_1\}} \frac{a_v + 1}{2} \operatorname{Var}_{\mu_{x,v,\gamma}}(\varphi) - \frac{1}{2} \operatorname{Var}_{\nu_{x,\gamma}}(\Delta).$$

PROOF. Let $x \in \Omega_{e_0}$. Using the definitions (5.17) and (5.18) of $x_v^{(\gamma\varphi)}$ and $Y_T^{(\gamma\varphi)}(x)$, respectively, we find that these two quantities are differentiable in $\gamma$ with the following derivatives:

(6.9) $$\frac{\partial}{\partial \gamma} x_v^{(\gamma\varphi)} = \sum_{\substack{e \in E: \\ e \ni v}} \varphi(e) e^{\gamma\varphi(e)} x_e \quad \text{and}$$

(6.10) $$\frac{\partial}{\partial \gamma} Y_T^{(\gamma\varphi)}(x) = \sum_{e \in T} \varphi(e) Y_T^{(\gamma\varphi)}(x) = \Delta(T) Y_T^{(\gamma\varphi)}(x).$$

Hence, it follows from the explicit form of $f_\gamma(x)$ given in (5.20) that $\gamma \mapsto f_\gamma(x)$ is continuously differentiable with

$$\frac{\partial f_\gamma}{\partial \gamma}(x) = \frac{a_{v_1}}{2} + \frac{1}{4} - \sum_{e \in E} \varphi(e) a_e + \sum_{v \in V \setminus \{v_0, v_1\}} \frac{a_v + 1}{2 x_v^{(\gamma\varphi)}} \sum_{\substack{e \in E: \\ e \ni v}} \varphi(e) e^{\gamma\varphi(e)} x_e$$
(6.11)
$$- \frac{1}{2} \frac{\sum_{T \in \mathcal{T}} \Delta(T) Y_T^{(\gamma\varphi)}(x)}{\sum_{T \in \mathcal{T}} Y_T^{(\gamma\varphi)}(x)}.$$

Inserting the definition of the measures $\mu_{x,v,\gamma}$ and $\nu_{x,\gamma}$ gives (6.7). Differentiating (6.11) with respect to $\gamma$ yields the claim for the second derivative. □

Now we are ready to prove the bound of the relative entropy:

PROOF OF LEMMA 6.1. By Assumption 2.3, $\varphi(e) \in [0, 1]$ for all $e \in E$ and, hence, $\Delta(T) \in [0, |E|]$ for all $T \in \mathcal{T}$. Thus, it follows from Lemma 6.2 that

(6.12) $$\sup_{x \in \Omega_{e_0}} \left| \frac{\partial^j f_\gamma}{\partial^j \gamma}(x) \right| < \infty$$

for $j = 1, 2$. Consequently, $g(\gamma) = E_{\mathbb{P}_{v_0,v_1,e_0}}[f_\gamma]$ is twice continuously differentiable and the derivatives are obtained as expectations of the derivatives of $f_\gamma$. Since $g(\gamma)$ is a relative entropy by its definition (6.2), one has $g(\gamma) \geq 0$ for all $\gamma$. Furthermore, $g(0) = 0$ because $\Pi^{(\gamma\varphi)}_{v_0,v_1,e_0} = \mathbb{P}_{v_0,v_1,e_0}$ for $\gamma = 0$. Thus, $g'(0) = 0$, and a Taylor expansion of $g$ around 0 yields

$$(6.13) \qquad g(\gamma) = \int_0^\gamma g''(\tilde\gamma)(\gamma - \tilde\gamma)\, d\tilde\gamma$$

with

$$(6.14) \qquad g''(\tilde\gamma) = E_{\mathbb{P}_{v_0,v_1,e_0}}\left[\frac{\partial^2 f_{\tilde\gamma}}{\partial^2 \tilde\gamma}\right].$$

To bound $g''$, we estimate the expression (6.8) for the second derivative of $f_\gamma$: Since variances are nonnegative, we obtain for all $\tilde\gamma \in \mathbb{R}$ the bound

$$(6.15) \qquad \begin{aligned} \frac{\partial^2 f_{\tilde\gamma}}{\partial^2 \tilde\gamma}(x) &\leq \sum_{v \in V} \frac{a_v + 1}{2} \operatorname{Var}_{\mu_{x,v,\tilde\gamma}}(\varphi) \\ &\leq \sum_{v \in V} \frac{a_v + 1}{2} \max_{\substack{e, e' \in E: \\ v \in e, v \in e'}} (\varphi(e) - \varphi(e'))^2 = S_\varphi. \end{aligned}$$

Since the upper bound $S_\varphi$ is independent of $x$ and $\tilde\gamma$, the claim of the lemma follows from (6.14) and (6.13). □

**7. Proof of the abstract version of the main lemma.** In this section, we combine the results from Sections 4–6 to prove the abstract version of the main lemma:

PROOF OF LEMMA 2.5. Combining (4.13) with the key estimate (5.16) yields

$$(7.1) \quad \log E_{\mathbb{Q}_{v_0,e_0}}\left[\left(\frac{x_{v_1}}{x_{v_0}}\right)^{1/4}\right] \leq E_{\Pi^{(\gamma\varphi)}_{v_0,v_1,e_0}}[H_{v_0,v_1}] + E_{\Pi^{(\gamma\varphi)}_{v_0,v_1,e_0}}\left[\log \frac{d\Pi^{(\gamma\varphi)}_{v_0,v_1,e_0}}{d\mathbb{P}_{v_0,v_1,e_0}}\right]$$

for all $\gamma \in \mathbb{R}$. To estimate the upper bound, we use Lemma 5.7 for the first term and the bound (6.1) of the relative entropy from Lemma 6.1 for the second term:

$$(7.2) \qquad \log E_{\mathbb{Q}_{v_0,e_0}}\left[\left(\frac{x_{v_1}}{x_{v_0}}\right)^{1/4}\right] \leq \frac{\gamma}{4} + \frac{S_\varphi \gamma^2}{2}.$$

This upper bound is minimal for $\gamma = -1/(4S_\varphi)$. Inserting this value for $\gamma$ gives the bound

$$(7.3) \qquad \log E_{\mathbb{Q}_{v_0,e_0}}\left[\left(\frac{x_{v_1}}{x_{v_0}}\right)^{1/4}\right] \leq -\frac{1}{32 S_\varphi}.$$

This proves the claim (2.19) of Lemma 2.5. □

**8. Application to finite boxes.** In this section, we derive the "concrete version" of the main lemma (Lemma 2.2) from the "abstract version" (Lemma 2.5). Given $\ell \in L_r$ (sufficiently far from the origin) and a sufficiently large box size $i \geq i_0(\ell)$ with $i_0(\ell)$ introduced immediately before Lemma 2.2, we start by defining a specific function $\varphi : E_r^{(i)} \to [0,1]$ as in Assumption 2.3(d) in the special case of finite boxes with periodic boundary conditions. Recall the definition of the "Dirichlet form" $S_\varphi$ in (2.20). We prove the following:

LEMMA 8.1. *Let $r \in \mathbb{N}$ with $r \geq 130$ and let the initial weight satisfy $a \in (0, (r-129)/256)$. There exist constants $l_0 = l_0(r,a) \in \mathbb{N}$ and $\xi = \xi(r,a) > 0$, such that for all $\ell \in L_r$ with $|\ell|_\infty \geq l_0$ and $i \geq i_0(\ell)$, there is a function $\varphi = \varphi_\ell^{(i)} : E_r^{(i)} \to [0,1]$ with the following properties:*

(a) *Assumption 2.3 is fulfilled for the graph $(V, E) = (V_r^{(i)}, E_r^{(i)})$ with constant initial weights $a$, the vertices $v_0 = \mathbf{0}$ and $v_1 = \ell$, the reference edge $e_0 = \{\mathbf{0}, (1,0)\}$, and the function $\varphi$.*
(b) *The following bound holds:*

$$(8.1) \qquad S_\varphi \leq \frac{1}{32(1+\xi)\log(|\ell|_\infty/r)}.$$

From now on, let $r \in \mathbb{N}$ with $r \geq 2$. Only later, in the proof of Lemma 8.1, we need the stronger bound $r \geq 130$.

*Defining a function $\varphi$.* For $u', v' \in L_r$ with $|u' - v'| = r$, we call $\{u', v'\}$ an *$r$-edge*. We say that an edge $e = \{u,v\} \in E_r$ belongs to an $r$-edge $\{u', v'\}$ if $e \subseteq [u', v']$, where $[u', v']$ denotes the straight line with end points $u'$ and $v'$. For $v \in V_r$, we call

$$(8.2) \qquad \mathrm{level}(v) := \lceil |v|_\infty / r \rceil$$

the *level* of the vertex, and we will also say that $v$ is at level $\mathrm{level}(v)$. The definition of levels is illustrated in Figure 3.

Before we define $\varphi$, we introduce an auxiliary function $D$:

DEFINITION 8.2 (Approximate Green's function). Given the two end points $u, v \in V_r$ of an edge $\{u,v\} \in E_r$, we take the unique vertices $u', v' \in L_r$ and the unique numbers $j_u, j_v \in \{0, 1, \ldots, r-1\}$ such that the following hold:

(a) $\{u', v'\}$ is an $r$-edge,
(b) $j_u + j_v = r - 1$ and
(c)
$$(8.3) \qquad u = \frac{j_u + 1}{r} u' + \frac{j_v}{r} v' \quad \text{and} \quad v = \frac{j_v + 1}{r} v' + \frac{j_u}{r} u'.$$

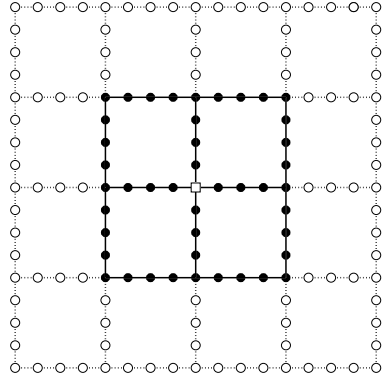

Fig. 3. *There is precisely one vertex at level* 0, *namely, the origin* **0**, *which is marked with a white square. The vertices at level* 1 *are marked with thick black dots, whereas the vertices at level* 2 *are marked with white circles.*

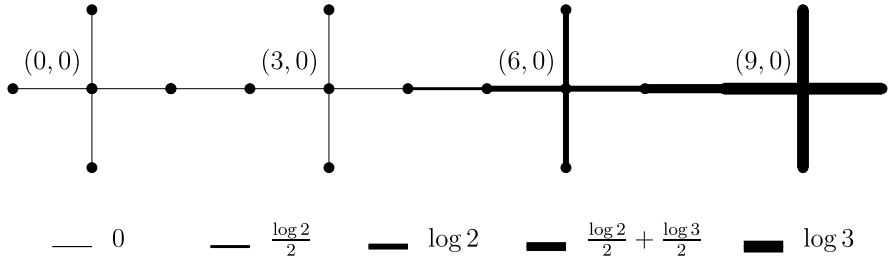

Fig. 4. *The values of the function $D$ on a piece of $G_3$.*

In particular, the straight line $[u', v']$ corresponding to the $r$-edge $\{u', v'\}$ contains the straight line $[u, v]$ corresponding to the edge $\{u, v\}$. We define

$$(8.4) \quad D(u,v) := \frac{j_u}{r-1} \log(\text{level}(u') \vee 1) + \frac{j_v}{r-1} \log(\text{level}(v') \vee 1).$$

Interchanging the roles of $u$ and $v$ interchanges $u'$ with $v'$ and $j_u$ with $j_v$. As a consequence, $D(u,v) = D(v,u)$. Thus, $D$ can be viewed as a function on *undirected* edges:

$$(8.5) \quad D: E_r \to [0, \infty), \quad D(\{u,v\}) := D(u,v).$$

REMARK 8.3.

(a) The map $D$ linearly interpolates on edges between two neighboring vertices in $L_r$; see Figure 4.
(b) For a technical reason we multiply the first [respectively second] log in (8.4) by $j_u/(r-1)$ [respectively $j_v/(r-1)$] and not by $(j_u+1/2)/r$ [respectively $(j_v+1/2)/r$]: This choice makes the function $D$ constant on all four edges incident to any vertex of degree 4. This feature is inherited by the function $\varphi$, defined below. It is essential that this property holds for the vertices $\mathbf{0}$ and $\ell$ by Assumption 2.3(d). It is not an essential property for all other vertices of degree 4, but nice to have, because it makes the terms corresponding to these vertices cancel out when estimating $S_\varphi$; see the arguments leading to (8.9) below.

For $v \in V_r$, define

(8.6) $$\underline{D}(v) := \min_{\substack{e \in E_r: \\ e \ni v}} D(e).$$

Note that $\underline{D}(\ell) > 0$ for all $\ell \in L_r$ with $\operatorname{level}(\ell) \geq 2$.

For fixed $\ell \in L_r$, we next rescale and truncate $D$ in such a way that it takes only values between 0 and 1, and the value 1 is taken on all edges which have both endpoints on a level $\geq \operatorname{level}(\ell)$.

We now use the identifications (2.3)–(2.5) of vertices/edges in finite boxes $G_r^{(i)}$ with vertices/edges in the infinite graph $G_r$. Recall, in particular, the definition (2.4) of $\hat{E}_r^{(i)}$.

DEFINITION 8.4. For $\ell \in L_r$ with $\underline{D}(\ell) > 0$, $i \geq i_0(\ell)$, and $e \in E_r^{(i)}$, define

(8.7) $$\varphi(e) := \varphi_\ell^{(i)}(e) := \frac{D(e)}{\underline{D}(\ell)} \wedge 1 \qquad \text{if } e \in \hat{E}_r^{(i)},$$

and $\varphi(e) := 1$ for periodically closing edges $e \in E_r^{(i)} \setminus \hat{E}_r^{(i)}$.

Here is the motivation why to use this particular function $\varphi$: Recall that $\varphi$ needs to be chosen in such a way that the "Dirichlet form" $S_\varphi$, defined in (2.20), gets small enough, as specified in (8.1). Given the constraints on $\varphi$ from Assumption 2.3(d), the Dirichlet form $S_\varphi$ is minimized for a certain two-point "Green's function." For our purposes, it suffices to use an approximation to this function having a roughly similar logarithmic decay. On the other hand, an arbitrary "Urysohn-type" function would in general not yield the required bound (8.1) of $S_\varphi$ roughly proportional to $1/\log |\ell|_\infty$.

PROOF OF LEMMA 8.1. We first verify Assumptions 2.3(a)–(d). Let $r \in \mathbb{N}$ with $r \geq 130$ and $a \in (0, (r-129)/256)$. Let $\ell \in L_r$ with $|\ell|_\infty \geq l_0$ and $i \geq i_0(\ell)$. Since $r \geq 2$, $\mathbf{0}$ and $\ell$ are not adjacent. Hence, Assumption 2.3(a) is satisfied. The reference edge $e_0 = \{\mathbf{0}, (1,0)\}$ is incident to $\mathbf{0}$; thus, 2.3(b) holds. Remark 2.4 shows that part (c) of the assumption holds. By definition, $\varphi$ takes only values in $[0, 1]$. To verify that part (d) of the assumption holds, note that each edge $e \in E_r^{(i)}$ incident to $\mathbf{0}$ belongs to an $r$-edge of the form $\{\mathbf{0}, \ell'\}$ with $\ell' \in L_r$ and $\text{level}(\ell') = 1$. Hence, for all $e \in E_r^{(i)}$ with $e \ni \mathbf{0}$, $D(e) = 0$ and $\varphi(e) = 0$. Furthermore, for $i \geq i_0(\ell)$, one has $D(e) \geq \underline{D}(\ell)$ for all $e \in E_r^{(i)}$ adjacent to $\ell$ and, thus, $\varphi(e) = 1$ for these edges. Thus, Assumption 2.3(d) holds. This completes the proof of part (a) of the lemma.

We now prove part (b) of the lemma. Let $d_v$ denote the vertex degree of a vertex $v$ in the graph $G_r^{(i)}$, that is, the number of edges incident to $v$ in $G_r^{(i)}$. Furthermore, let

$$\alpha := \alpha(a) := \frac{2a+1}{2}. \tag{8.8}$$

In the special case of the graph $G_r^{(i)}$ and with our choice of $\varphi$, for all vertices $v$ of vertex degree $d_v = 4$, that is, $v \in L_r^{(i)}$, we have $\varphi(e) - \varphi(e') = 0$ for all edges $e, e' \in E_r^{(i)}$ incident to $v$. Furthermore, all vertices in $V_r^{(i)}$ have either degree 2 or 4. Hence, we get

$$\begin{aligned}S_\varphi &= \sum_{v \in V_r^{(i)}} \frac{a d_v + 1}{2} \max_{\substack{e,e' \in E_r^{(i)}: \\ v \in e, v \in e'}} (\varphi(e) - \varphi(e'))^2 \\ &= \alpha \sum_{v \in V_r^{(i)} \setminus L_r^{(i)}} \max_{\substack{e,e' \in E_r^{(i)}: \\ v \in e, v \in e'}} (\varphi(e) - \varphi(e'))^2.\end{aligned} \tag{8.9}$$

By Remark 8.3, the function $\tilde{E}_r^{(i)} \ni e \mapsto D(e)$ interpolates linearly between the values $\log(\text{level}(u') \vee 1)$ and $\log(\text{level}(v') \vee 1)$ for edges $e$ contained in the convex hull of an $r$-edge $\{u', v'\}$. As a consequence, for $v \in V_r^{(i)}$ with $\text{level}(v) \geq \text{level}(\ell) + 1$ and all edges $e \ni v$, one has $D(e) \geq \underline{D}(\ell)$; recall the definition (8.6) of $\underline{D}$. This implies $\varphi(e) = 1$. Thus, we obtain

$$\begin{aligned}S_\varphi &\leq \alpha \sum_{\substack{v \in V_r^{(i)} \setminus L_r^{(i)}: \\ \text{level}(v) \leq \text{level}(\ell)}} \max_{\substack{e \in E_r^{(i)}: \\ v \in e}} \left(\varphi(e) - \min_{\substack{e' \in E_r^{(i)}: \\ v \in e'}} \varphi(e')\right)^2 \\ &= \alpha \sum_{\substack{v \in V_r^{(i)} \setminus L_r^{(i)}: \\ \text{level}(v) \leq \text{level}(\ell)}} \max_{\substack{e \in E_r^{(i)}: \\ v \in e}} \left(\frac{D(e)}{\underline{D}(\ell)} - \min_{\substack{e' \in E_r^{(i)}: \\ v \in e'}} \frac{D(e')}{\underline{D}(\ell)}\right)^2 \\ &= \frac{\alpha}{\underline{D}(\ell)^2} \sum_{\substack{v \in V_r^{(i)} \setminus L_r^{(i)}: \\ \text{level}(v) \leq \text{level}(\ell)}} \max_{\substack{e \in E_r^{(i)}: \\ v \in e}} (D(e) - \underline{D}(v))^2.\end{aligned} \tag{8.10}$$

If $e = \{u, v\} \in E_r$ is adjacent to a vertex $\ell' \in L_r$, say, $v = \ell'$, then in the representation (8.3) for $u$ and $v$, one has $j_u = 0$ and $j_v = r - 1$. Hence, for all $\ell' \in L_r$,

(8.11) $$\underline{D}(\ell') = \log(\text{level}(\ell') \vee 1).$$

Assume that $v \in V_r^{(i)} \setminus L_r^{(i)}$. Then, there exist precisely two edges $e = \{u, v\}$ and $\tilde{e} = \{\tilde{u}, v\}$ incident to $v$, and both of them belong to the same $r$-edge. Let $u', v' \in L_r$ denote the endpoints of this $r$-edge as in (8.3). The definition (8.4) and (8.5) of $D$ and (8.11) imply

(8.12) $$\max_{\substack{e \in E_r^{(i)}: \\ e \ni v}} (D(e) - \underline{D}(v))^2 = \left(\frac{\underline{D}(u') - \underline{D}(v')}{r-1}\right)^2.$$

The convex hull $[u', v']$ of any $r$-edge $e' = \{u', v'\}$ contains precisely $r+1$ vertices in $V_r$, where two among them belong to $L_r$. We get

(8.13) $$\sum_{v \in (V_r^{(i)} \setminus L_r^{(i)}) \cap [u', v']} \max_{\substack{e \in E_r^{(i)}: \\ e \ni v}} (D(e) - \underline{D}(v))^2 = (r-1)\left(\frac{\underline{D}(u') - \underline{D}(v')}{r-1}\right)^2$$
$$= \frac{(\underline{D}(u') - \underline{D}(v'))^2}{r-1}.$$

This formula is the essential point where one can see why large values of $r$ make recurrence easier to prove: Every $r$-edge contributes $r-1$ summands; each of them scales proportional to $1/(r-1)^2$. Thus, every $r$-edge gives a contribution proportional to $1/(r-1)$, which gets small for large values of $r$.

Inserting the expression (8.13) into (8.10), we obtain

(8.14) $$S_\varphi \leq \frac{\alpha}{\underline{D}(\ell)^2} \sum_{\substack{\{u', v'\} r\text{-edge:} \\ \text{level}(u'), \text{level}(v') \leq \text{level}(\ell)}} \frac{(\underline{D}(u') - \underline{D}(v'))^2}{r-1}.$$

Before going into the details, we explain the rough idea of the remainder of the proof: The summand corresponding to an $r$-edge $\{u', v'\}$ in the sum in (8.14) is roughly of the order $l^{-2}(r-1)^{-1}$ with $l = \text{level}(u')$. The number of $r$-edges on this level has roughly the size $\text{const} \cdot l$. Thus, the contribution from edges on level $l$ in the sum in (8.14) is roughly $\text{const} \cdot l^{-1}(r-1)^{-1}$. Thus, the whole sum [without the first factor $\alpha/\underline{D}(\ell)^2$] has roughly the order of magnitude $\text{const} \cdot \log(\text{level}(\ell))(r-1)^{-1}$. The first factor in (8.14) equals $\alpha/(\log \text{level}(\ell))^2$. Hence, we get an upper bound for $S_\varphi$ roughly of the size $\text{const}/[(r-1)\log \text{level}(\ell)]$.

More precisely, we proceed as follows: Note that $r$-edges $\{u', v'\}$ with $\text{level}(u') = \text{level}(v')$ do not contribute to the sum in (8.14). Furthermore, there are $4(2l - 1)$ $r$-edges connecting points $u' \in L_r$ with $\text{level}(u') = l$ to points $v' \in L_r$ with $\text{level}(v') = l - 1$. Using (8.11), we obtain

$$
\begin{aligned}
S_\varphi &\leq \frac{\alpha}{(r-1)[\log(\text{level}(\ell))]^2} \\
&\quad \times \sum_{l=1}^{\text{level}(\ell)} 4(2l-1)[\log(l \vee 1) - \log((l-1) \vee 1)]^2 \\
&= \frac{4\alpha}{(r-1)[\log(\text{level}(\ell))]^2} \sum_{l=2}^{\text{level}(\ell)} (2l-1)\left[\log \frac{l}{l-1}\right]^2.
\end{aligned}
\tag{8.15}
$$

Applying the inequality $\log x \leq x - 1$ yields

$$
\begin{aligned}
\sum_{l=2}^{\text{level}(\ell)} (2l-1)\left[\log \frac{l}{l-1}\right]^2 &\leq \sum_{l=2}^{\text{level}(\ell)} (2l-1)\left[\frac{l}{l-1} - 1\right]^2 \\
&= 2 \sum_{l=1}^{\text{level}(\ell)-1} \frac{1}{l} + \sum_{l=1}^{\text{level}(\ell)-1} \frac{1}{l^2} \\
&\leq 2[1 + \log(\text{level}(\ell))] + 2.
\end{aligned}
\tag{8.16}
$$

Combining this with (8.15), it follows that

$$
S_\varphi \leq \frac{8\alpha[\log(\text{level}(\ell)) + 2]}{(r-1)[\log(\text{level}(\ell))]^2}.
\tag{8.17}
$$

Using the definition (8.8) of $\alpha$ and our assumption on the initial weight $a$, we get

$$
\frac{r-1}{256\alpha} = \frac{r-1}{128(2a+1)} > 1.
\tag{8.18}
$$

Consequently, we can choose a constant $c = c(r, a) \in (0, 1)$ such that

$$
\frac{c(r-1)}{256\alpha} > 1
\tag{8.19}
$$

holds. We set

$$
\xi = \xi(r, a) := \frac{c(r-1)}{256\alpha} - 1.
\tag{8.20}
$$

By (8.19), $\xi > 0$. We choose $l_0(c) = l_0(r, a)$ large enough that for all $\ell \in L_r$ with $|\ell|_\infty \geq l_0(c)$, one has

$$
\frac{\log(\text{level}(\ell)) + 2}{\log(\text{level}(\ell))} \leq c^{-1}.
\tag{8.21}
$$

Inserting this bound into (8.17) and using the definition of $\xi$, we obtain

$$(8.22) \qquad S_\varphi \le 1\frac{8\alpha}{c(r-1)\log(\text{level}(\ell))} = \frac{1}{32(1+\xi)\log(\text{level}(\ell))}.$$

Using the relation $\text{level}(\ell) = |\ell|_\infty/r$ yields the claim of the lemma. $\square$

Finally, we deduce the concrete version of the main lemma:

PROOF OF LEMMA 2.2. Let $r \in \mathbb{N}$ with $r \ge 130$ and $a \in (0, (r-129)/256)$. By Lemmas 8.1 and 2.5, there exist constants $l_0 = l_0(r,a) \in \mathbb{N}$ and $\xi = \xi(r,a) > 0$, such that for all $\ell \in L_r$ with $|\ell|_\infty \ge l_0$ and $i \ge i_0(\ell)$, one has

$$(8.23) \qquad E_{\mathbb{Q}_\mathbf{0}^{(i)}}\left[\left(\frac{x_\ell}{x_\mathbf{0}}\right)^{1/4}\right] \le \exp\left(-\frac{1}{32 S_\varphi}\right) \le \left(\frac{r}{|\ell|_\infty}\right)^{1+\xi}. \qquad \square$$

**9. Proving the main theorems.** Now we are ready to prove the estimates for the hitting probabilities of the edge-reinforced random walk.

PROOF OF THEOREM 1.2. Given $r$ and $a$ as in the hypothesis of the theorem, take $l_0 = l_0(r, a)$ and $\xi = \xi(r, a)$ from Lemma 2.2. Let $\ell \in L_r$ with $|\ell|_\infty \ge l_0$. Fix any finite set $\Pi$ of finite paths in $G_r$ starting in $\mathbf{0}$ and ending in $\ell$, such that every path $\pi \in \Pi$ visits $\mathbf{0}$ only once, namely, at its start, and visits also $\ell$ only once, namely, at its end. Take any $i \ge i_0(\ell)$ large enough so that for all $\pi \in \Pi$ and for all vertices $v$ in $\pi$, the vertex $v$ and all its neighboring vertices belong to $\tilde{V}_r^{(i)}$. Identifying $\tilde{V}_r^{(i)}$ with a subset of $V_r^{(i)}$, we view $\Pi$ also as a set of finite paths in $G_r^{(i)}$.

For any admissible finite path $\pi = (v_0, v_1, \ldots, v_n)$ in $G_r^{(i)}$, let $\pi^\leftrightarrow = (v_n, \ldots, v_1, v_0)$ denote its reversed path. We set $\Pi^\leftrightarrow = \{\pi^\leftrightarrow : \pi \in \Pi\}$. Furthermore, let $A_\pi = \{X_t = v_t \text{ for } t = 0, \ldots, n\}$ denote the event that the random walk $(X_t)_t$ follows the path $\pi$ as an initial piece of its trajectory. By a slight abuse of notation, we use this notation both for random walks on $G_r$ and random walks on $G_r^{(i)}$. Set $A_\Pi = \bigcup_{\pi \in \Pi} A_\pi$ and, similarly, $A_{\Pi^\leftrightarrow} = \bigcup_{\pi \in \Pi^\leftrightarrow} A_\pi$.

Fix a collection of weights $x \in (0, \infty)^{E_r^{(i)}}$. Let $Q_{v_0, x}^{(i)}$ denote the law of the Markovian random walk on $G_r^{(i)}$ starting in $v_0 \in V_r^{(i)}$ as defined in (2.11) and (2.12). We are going to compare the laws $Q_{\mathbf{0},x}^{(i)}$ and $Q_{\ell,x}^{(i)}$ of the two random walks starting in $\mathbf{0}$ and $\ell$, respectively: For a finite path $\pi = (v_0, v_1, \ldots, v_n)$ from $v_0 = \mathbf{0}$ to $v_n = \ell$ as above, we have

$$(9.1) \quad Q_{\mathbf{0},x}^{(i)}[A_\pi] = \prod_{t=0}^{n-1} \frac{x_{\{v_t, v_{t+1}\}}}{x_{v_t}} = \frac{x_\ell}{x_\mathbf{0}} \prod_{t=1}^{n} \frac{x_{\{v_t, v_{t-1}\}}}{x_{v_t}} = \frac{x_\ell}{x_\mathbf{0}} Q_{\ell,x}^{(i)}[A_{\pi^\leftrightarrow}]$$

and, thus,

$$(9.2) \qquad Q_{\mathbf{0},x}^{(i)}[A_\Pi] = \frac{x_\ell}{x_\mathbf{0}} Q_{\ell,x}^{(i)}[A_{\Pi^\leftrightarrow}] \le \frac{x_\ell}{x_\mathbf{0}}.$$

Trivially, we know $Q_{\mathbf{0},x}^{(i)}[A_\Pi] \leq 1$ and, thus,

$$(9.3) \qquad Q_{\mathbf{0},x}^{(i)}[A_\Pi] \leq \frac{x_\ell}{x_\mathbf{0}} \wedge 1 \leq \left(\frac{x_\ell}{x_\mathbf{0}}\right)^{1/4}.$$

Representing the reinforced random walk on $G_r^{(i)}$ with start in $\mathbf{0}$ as a mixture of the Markovian random walks $Q_{\mathbf{0},x}^{(i)}$ (Lemma 2.1), we conclude, using the claim (2.16) of Lemma 2.2 in the last step,

$$(9.4) \quad P_{\mathbf{0},a}^{(i)}[A_\Pi] = E_{\mathbb{Q}_\mathbf{0}^{(i)}}[Q_{\mathbf{0},x}^{(i)}[A_\Pi]] \leq E_{\mathbb{Q}_\mathbf{0}^{(i)}}\left[\left(\frac{x_\ell}{x_\mathbf{0}}\right)^{1/4}\right] \leq \left(\frac{r}{|\ell|_\infty}\right)^{1+\xi}.$$

Now, the reinforced random walk in $G_r$ and the reinforced random walk in $G_r^{(i)}$ with the same starting point $\mathbf{0}$ and the same initial weights $a$ have the same law up to the stopping time when the random walker leaves the graph $(\tilde{V}_r^{(i)}, \tilde{E}_r^{(i)})$. Since $i$ was chosen large enough, depending on $\Pi$, this implies

$$(9.5) \qquad P_{\mathbf{0},a}^{G_r}[A_\Pi] = P_{\mathbf{0},a}^{(i)}[A_\Pi] \leq \left(\frac{r}{|\ell|_\infty}\right)^{1+\xi}.$$

Now let $\Pi_{\mathbf{0},\ell}$ denote the set of *all* finite paths in $G_r$ from $\mathbf{0}$ to $\ell$ that visit $\mathbf{0}$ and $\ell$ only once. We write $\Pi_{0,\ell}$ as the union of an increasing sequence of finite subsets $\Pi_n \uparrow \Pi_{0,\ell}$, $n \to \infty$. Then as $n \to \infty$, the events $A_{\Pi_n}$ also increase toward the event $\{\tau_\ell < \tau_\mathbf{0}\}$ that the random walker hits $\ell$ before returning to $\mathbf{0}$. We get

$$(9.6) \qquad P_{\mathbf{0},a}^{G_r}[\tau_\ell < \tau_\mathbf{0}] = \lim_{n \to \infty} P_{\mathbf{0},a}^{G_r}[A_{\Pi_n}] \leq \left(\frac{r}{|\ell|_\infty}\right)^{1+\xi},$$

which proves the first part of the theorem.

To prove the second part, for $l \in \mathbb{N}$, let $\mathcal{L}_l$ denote the set of vertices $\ell \in L_r$ with $|\ell|_\infty = rl$. Observe that the first time the reinforced random walker starting at $\mathbf{0}$ hits the set $\mathcal{V}_l$ (provided this ever happens), it hits it in a vertex in $\mathcal{L}_l$. This implies the following bound for the probability that the edge-reinforced random walk visits the set $\mathcal{V}_l$ before returning to $\mathbf{0}$:

$$(9.7) \qquad P_{\mathbf{0},a}^{G_r}[\tau_{\mathcal{V}_l} < \tau_\mathbf{0}] = P_{\mathbf{0},a}^{G_r}[\tau_{\mathcal{L}_l} < \tau_\mathbf{0}] \leq \sum_{\ell \in \mathcal{L}_l} P_{\mathbf{0},a}^{G_r}[\tau_\ell < \tau_\mathbf{0}].$$

Since there are $8l$ vertices in $L_r$ with $|\ell|_\infty = rl$ and hence in $\mathcal{L}_l$, it follows from the first part of the theorem that

$$(9.8) \qquad P_{\mathbf{0},a}^{G_r}[\tau_{\mathcal{V}_l} < \tau_\mathbf{0}] \leq 8l\left(\frac{1}{l}\right)^{1+\xi} = 8l^{-\xi}.$$

This concludes the proof of the theorem. $\square$

Our estimates for the hitting probabilities imply recurrence:

PROOF OF THEOREM 1.1. The edge-reinforced random walk on any infinite, locally finite, connected graph visits infinitely many points with probability one. In other words, it does not get "stuck" on a proper finite subgraph. This fact is, for example, remarked in [5], page 2. Here is the intuitive argument: Take two neighboring vertices $v$ and $v'$. Conditional on the random walk visiting $v$ the $k$th time at time $t$, the probability to be at $v'$ at time $t+1$ is at least $a_{\{v,v'\}}/(a_v + 2k)$, which is not summable in $k$. Thus, the following implication holds almost surely by a Borel–Cantelli argument: if $v$ is visited infinitely often, then $v'$ is also visited infinitely often.

Hence, the probability that the edge-reinforced random walk never returns to $\mathbf{0}$ equals the probability that the random walker visits the sets $\mathcal{V}_l$ for all $l$ before returning to $\mathbf{0}$. This shows that

$$(9.9) \qquad P^{G_r}_{\mathbf{0},a}[\tau_{\mathbf{0}} = \infty] = P^{G_r}_{\mathbf{0},a}\left[\bigcap_{l=1}^{\infty} \{\tau_{\mathcal{V}_l} < \tau_{\mathbf{0}}\}\right] = \lim_{l \to \infty} P^{G_r}_{\mathbf{0},a}[\tau_{\mathcal{V}_l} < \tau_{\mathbf{0}}].$$

By estimate (1.10) from Theorem 1.2, the last limit equals zero. Hence, $P^{G_r}_{\mathbf{0},a}[\tau_{\mathbf{0}} = \infty] = 0$, which means that the edge-reinforced random walk returns to its starting point with probability one. By Theorem 2.1 of [8], this implies that the edge-reinforced random walk visits every vertex infinitely often with probability one. □

**10. Conclusion.** The technique used in this paper is quite robust to some perturbations of the model, but very nonrobust to other perturbations. More specifically:

- The method of the proof in this paper could be adapted to the $r$-diluted versions of a large class of two-dimensional graphs with constant or periodic initial weights, provided $r$ is large. The graphs need to have sufficiently many reflection symmetries, as specified in Assumption 2.3. This includes, among others, the diluted versions of $\mathbb{Z}^2$ with additional diagonal edges added in a periodic and reflection-symmetric way. It also includes diluted versions of $\mathbb{Z}^2 \times G$, with any finite graph $G$ with a transitive action of the automorphism group, for example, $\mathbb{Z}^2 \times (\mathbb{Z} \bmod N\mathbb{Z})$. The abstract version of our main Lemma 2.5 is general enough to cover all these cases.
- A variant of the method presented here was used in [7] to analyze the edge-reinforced random walk on a large class of one-dimensional graphs with a reflection symmetry. There, besides recurrence, exponential decay of the weights of the random environment was proven. The abstract version of our main Lemma 2.5 also covers these cases. In the one-dimensional setup, the law of the position of the reinforced random walk is tight uniformly in time with exponential tails. The proof of this fact is contained in [7]. Indeed, in one dimension, the 1/4-moment analog to (2.16) decays not only with a power law, but even exponentially fast. Furthermore, such a bound holds also for the one-dimensional infinite-volume random environment. In the present paper, we do not treat a two-dimensional analog.
- On the other hand, we cannot handle arbitrary initial weights with our method, as soon as one violates reflection symmetry.
- The techniques of this paper break immediately down if one varies the reinforcement scheme. They apply only to linear reinforcement. The method depends crucially on partial exchangeability of the linearly edge-reinforced random walk and, even more, on the specific form of the random environment distribution as described in Lemma 3.4 and Definition 3.2.
- The method presented here applies also to the original problem on $\mathbb{Z}^2$ and yields bounds for the random environment distribution because the abstract version of the main lemma (Lemma 2.5) applies in this case as well. One could use this to derive bounds for the hitting probabilities in the style of Theorem 1.2, but weaker. However, using the variational principle in a different way [using $\mathbb{Q}_{v_1,e_0}$ instead of $\mathbb{P}_{v_0,v_1,e_0}$ in the definition (5.3) of $F$] and using more sophisticated deformation maps (defined in terms of a function $\varphi$ as in Assumption 2.3, but additionally depending on the environment $x$), one can derive different and stronger bounds for the edge-reinforced random walk on $\mathbb{Z}^2$ beyond the scope of this paper. This includes a fast decay of the expected *logarithm* of the weights of the random environment for small initial weights $a$. For more details, see [9], in particular, Theorem 2.3 in this reference.

An important difference between the present paper and the previous work of [9] and [7] consists in the simplifying idea to prove recurrence directly from bounds on the environment on finite boxes without studying the random environment in infinite volume.

**Acknowledgments.** S. R. would like to thank Wolfgang Woess for very useful discussions during her visit to Graz. The authors would like to thank two anonymous referees for useful suggestions improving the presentation of this paper.

Mathematical Institute
University of Munich
Theresienstr. 39, D-80333 Munich
Germany
E-mail: merkl@mathematik.uni-muenchen.de

Zentrum Mathematik
Bereich M5
Technische Universität München
D-85747 Garching bei München
Germany
E-mail: srolles@ma.tum.de